\documentclass[12pt]{article}
\usepackage{amsmath,amssymb}
\usepackage{graphicx}
\usepackage{color}
\newcommand{\eqdef}{\stackrel{\text{def}}{=}}

\newcommand{\n}{\nonumber\\}
\newcommand{\bm}{\boldsymbol}
\newcommand{\ignore}[1]{}
\numberwithin{equation}{section}
\newcommand{\ai}{\text{I}}
\newcommand{\ait}{\text{II}}
\newcommand{\at}{\text{III}}
\newcommand{\cF}{c_{\text{\tiny$\mathcal{F}$}}}
\newcommand{\Romannumeral}[1]{\uppercase\expandafter{\romannumeral#1}}
\newcommand{\I}{\text{\Romannumeral{1}}}
\newcommand{\II}{\text{\Romannumeral{2}}}
\newtheorem{theo}{\bf Theorem}[section]

\newtheorem{coro}[theo]{\bf Corollary}

\newcommand{\ma}{\hspace{0pt}}
\allowdisplaybreaks[4]

\begin{document}

\baselineskip=20pt

\newcommand{\preprint}{
\vspace*{-20mm}
\begin{flushleft} Version:  October 28, 2019 \end{flushleft}
}
\newcommand{\Title}[1]{{\baselineskip=26pt
   \begin{center} \Large \bf #1 \\ \ \\ \end{center}}}
\newcommand{\Author}{\begin{center}
   \large \bf  Choon-Lin Ho${}^a$ and 
   Ryu Sasaki${}^{b}$ 
   \end{center}}
\newcommand{\Address}{\begin{center}
     $^a$  Department of Physics, Tamkang University, 
     Tamsui  25137, Taiwan (R.O.C.) \\
     ${}^b$Department of Physics, Tokyo University of Science,
     Noda 278-8510, Japan\\
    \end{center}}
\newcommand{\Accepted}[1]{\begin{center}
   {\large \sf #1}\\ \vspace{1mm}{\small \sf Accepted for Publication}
   \end{center}}

\thispagestyle{empty}

\Title{Discrete orthogonality relations for  multi-indexed Laguerre and Jacobi polynomials}                     %

\Author

\Address

\begin{abstract}
The discrete orthogonality relations hold for all the orthogonal polynomials 
obeying three term recurrence relations. 
We show that they also hold for multi-indexed Laguerre and Jacobi polynomials, 
which are new orthogonal polynomials obtained 
by deforming these classical orthogonal polynomials.
The discrete orthogonality relations could be considered as more encompassing characterisation 
of orthogonal polynomials than the three term recurrence relations.
As the multi-indexed orthogonal polynomials start at a positive degree $\ell_{\mathcal D}\ge1$,
the three term recurrence relations are broken.
The extra $\ell_{\mathcal D}$  `lower degree polynomials', 
which are necessary for the discrete orthogonality relations, are identified. 
The corresponding Christoffel
numbers are determined.
The main results are obtained by the blow-up analysis of the second order differential operators
governing the multi-indexed orthogonal polynomials around the zeros of these polynomials at a 
degree $\ell_{\mathcal D}+\mathcal{N}$. 
The discrete orthogonality relations are shown to hold for another group of `new' orthogonal
polynomials called Krein-Adler polynomials based on the Hermite, Laguerre and Jacobi polynomials.
\end{abstract}

\section{Introduction}
\label{intro}

The discrete orthogonality relations among orthogonal polynomials $\{P_n(\eta)\}$, 
satisfying the three term recurrence relations,
\begin{align*} 
\eta P_n(\eta)=a_n P_{n+1}(\eta)+b_n P_n(\eta)+c_n P_{n-1}(\eta),\quad n\ge0\quad
     \Leftrightarrow \int \!P_n(\eta)P_m(\eta)d\alpha(\eta)\propto \delta_{n\,m},
\end{align*}
 are formulated as follows. 
Here $d\alpha(\eta)$ is the orthogonality measure, which can be continuous or purely discrete.
Let us fix a positive integer $\mathcal{N}$ and denote the zeros of $P_{\mathcal N}(\eta)$ as
$\eta_1,\eta_2,\ldots,\eta_{\mathcal N}$, which are all distinct.
The discrete orthogonality relations hold among all the {\em lower degree polynomials\/}:
\begin{equation}
\sum_{j=1}^{\mathcal N}\lambda_jP_n(\eta_j)P_m(\eta_j)=0,\quad
n\neq m=0,1,\ldots,\mathcal{N}-1,
\label{disor}
\end{equation}
in which the {\em Christoffel numbers} $\{\lambda_j\}$  are defined as follows:
\begin{equation}
\lambda_j\eqdef \int \!\ell_j(\eta)d\alpha(\eta)>0,\quad
\ell_j(\eta)\eqdef \frac{P_{\mathcal N}(\eta)}{P'_{\mathcal N}(\eta_j)(\eta-\eta_j)},\quad
j=1,\ldots,\mathcal{N}.
\label{chris}
\end{equation}
The proof is very simple.  
It is well known that the Gaussian quadrature, 
\begin{equation*}
\int\!f(\eta)d\alpha(\eta)=\sum_{j=1}^{\mathcal N}\lambda_jf(\eta_j).
\end{equation*}
is exact for any polynomial $f(\eta)$
of degree at most $2\mathcal{N}-1$ (see, for example, (5.3.7) in \cite{askey}).
By putting $f(\eta)=P_n(\eta)P_m(\eta)$ and using the orthogonality of $\{P_n(\eta)\}$, 
we obtain \eqref{disor}.
 
 Recently, a new type of orthogonal polynomials, called {\em exceptional\/} and/or 
 {\em multi-indexed\/} orthogonal
 polynomials \cite{xop}-\cite{dismul}, 
 has been discussed by many researchers within the circles of ordinary differential (difference)
 equations, orthogonal polynomials and exactly solvable quantum mechanics.
 They are deformations of the classical orthogonal polynomials, 
 satisfying second order differential (difference) 
 equations but not the three term recurrence relations.
 As the classical orthogonal polynomials, 
 these new orthogonal polynomials span a complete orthogonal basis 
 in an appropriate Hilbert space.
 Typically  these polynomials, denoted as $\{P_{\mathcal{D},\,n}(\eta)\}$, $n=0,1,\ldots$, with
 $\mathcal{D}$ being the {\em multi-index}, see \eqref{defD},
 have   a non-vanishing lowest degree (see \eqref{ellD} for $\ell_{\mathcal D}$) 
 \begin{equation}
\text{deg}\bigl(P_{\mathcal{D},\,0}(\eta)\bigr)=\ell_{\mathcal D}\ge1,\quad
\text{deg}\bigl(P_{\mathcal{D},\,n}(\eta)\bigr)=\ell_{\mathcal D}+n,\quad n=1,2,\ldots.
\label{degrel}
\end{equation}
This breaks the three term recurrence relations and  it is essential for circumventing the constraints 
due to Routh-Bochner \cite{routh}-\cite{bochner}.

\bigskip
We ask the following questions.\\
{\bf Q1.} Is there a multi-indexed orthogonal polynomial version 
of the discrete orthogonality relations?\\
As above, by  fixing  a positive integer $\mathcal{N}$,  
$P_{\mathcal{D},\,\mathcal{N}}(\eta)$ has 
$\widetilde{\mathcal{N}}\eqdef \ell_{\mathcal D}+\mathcal{N}$ zeros, 
whereas there are only $\mathcal{N}$ ordinary 
lower degree polynomials $\{P_{\mathcal{D},\,n}(\eta)\}$, $n=0,1,\ldots, \mathcal{N}-1$.\\
{\bf Q2.} Is it possible to identify all the extra $\ell_{\mathcal D}$ `lower degree polynomials' 
which enter in the discrete orthogonality relations?\\
{\bf Q3.} Is it possible to determine the explicit forms of the Christoffel numbers $\{\lambda_j\}$ \eqref{chris}?

\bigskip
In the subsequent sections we will provide affirmative answers to all the above questions
by taking the best known examples of the new orthogonal polynomials, {\em i.e.\/}
the multi-indexed Laguerre and Jacobi polynomials \cite{os25,os37}
and the Krein-Adler (K-A) polynomials \cite{krein, adler} based 
on the Hermite, Laguerre and Jacobi polynomials.
Similar results are expected for the other new orthogonal polynomials, {\em e.g.} 
the multi-indexed Wilson, Askey-Wilson, Racah and $q$-Racah polynomials \cite{dismul} 
and the K-A polynomials based on these classical orthogonal polynomials, etc.
It is a good challenge to demonstrate them concretely.

A clue to answering these questions can be found in a recent paper \cite{zero}, 
which  will be referred to as the {\em zero perturbation\/} paper.
It investigated
the behaviour of the solutions of the Schr\"odinger-like equations containing  classical
orthogonal polynomials $\{P_n(\eta)\}$.
To be more precise, it  focused on the behaviour of the solutions around the zeros 
$\eta_1, \eta_2,\ldots,\eta_{\mathcal N}$ of $P_{\mathcal N}(\eta)$,
which was formulated as an eigenvalue problem of an 
$\mathcal{N}\times \mathcal{N}$ hermitian (real symmetric) matrix 
$\mathcal{M}$ obtained through perturbation.
The discrete orthogonality relations \eqref{disor} for the classical orthogonal polynomials were re-derived
as the orthogonality of the eigenvectors of $\mathcal{M}$ together with the explicit forms of the
Christoffel numbers.
As the multi-indexed orthogonal polynomials are the deformations of the classical orthogonal polynomials,
the reformulation of the zero perturbation paper in terms of the multi-indexed orthogonal polynomials
goes almost parallel. 
The above problems {\bf Q1}--{\bf Q3} are now cast in the framework of linear algebra
and the discrete orthogonality relations are reduced to the orthogonality of eigenvectors of a certain 
$\widetilde{\mathcal{N}}\times \widetilde{\mathcal{N}}$ matrix $\mathcal{M}$,
 which is {\em complex symmetric\/} rather than hermitian (real symmetric) 
as some of zeros $\eta_1, \eta_2,\ldots,\eta_{\widetilde{\mathcal N}}$ 
are now complex conjugate pairs.
By linear algebraic means, the discrete orthogonality ({\bf Q1}) is proved exactly and 
the extra `lower degree polynomials' ({\bf Q2}) and the Christoffel numbers ({\bf Q3}) are
determined by the basic equations 
\eqref{basic}, \eqref{rbasic}, \eqref{KAbasic}, \eqref{KArbasic},   
governing these polynomials.

\bigskip
This paper is organised as follows.
In section \ref{sec:mulLJ} the essence of the multi-indexed Laguerre and 
Jacobi polynomials are recapitulated 
in order to introduce appropriate notion and notation.
They are essentially a concise summary of  the original paper \cite{os25}. 
For the actual numerical verifications,
the newly formulated simplified expressions of the multi-indexed polynomials 
 in \cite{os37} are recommended.
Section \ref{sec:disortmulLJ} provides the main results of the paper. 
All three questions {\bf Q1}--{\bf Q3} are
answered affirmatively;
the explicit  forms of the extra $\ell_{\mathcal D}$ `lower degree polynomials' are presented
and their discrete orthogonality relations with the explicit forms of the Christoffel numbers
are demonstrated.
The concept of the extra `lower degree polynomials' is very clear and natural 
within the framework of the multi-indexed orthogonal polynomials.
The ordinary lower degree polynomials are obtained by reducing $\mathcal{N}$ to a lower degree $m$,
$P_{\mathcal{D},\,\mathcal{N}}(\eta)\to P_{\mathcal{D},\,m}(\eta)$, $m=0,\ldots,\mathcal{N}-1$,
whereas the extra `lower degree polynomials' are obtained by reducing one of the multi-index $d_j$ in 
$\mathcal{D}$ to an allowed lower degree $\epsilon_k$, 
$P_{\mathcal{D},\,\mathcal{N}}(\eta)\to P_{\mathcal{D}',\,\mathcal{N}}(\eta)$, 
see \eqref{edef}, \eqref{1D}.
The main results are presented as {\bf Theorem \ref{theo1}}. 
Section \ref{sec:zerosum}  provides the details of the tools to derive the main results 
presented in section \ref{sec:disortmulLJ}.
They are the eigenvalue problems of two
 $\widetilde{\mathcal{N}}\times \widetilde{\mathcal{N}}$
matrices $\widetilde{\mathcal{M}}$ and $\mathcal{M}$,  
and the basic equations satisfied by the extra `lower degree polynomials'. 
The starting point is the time-dependent Schr\"odinger equations of quantum mechanical systems 
corresponding to the multi-indexed Laguerre and Jacobi polynomials, 
a remake of the zero perturbation paper \cite{zero}.
The orthogonality of  the eigenvectors of $\mathcal{M}$ implies the discrete orthogonality, 
stated as {\bf Theorem \ref{theo:1}} and {\bf \ref{theo:2}}.
The new features, {\em e.g.\/} the complex conjugate zeros and non-hermitian inner products, etc are
emphasised.
The basic equations \eqref{basic} governing the behaviours of the extra `lower degree polynomials'
are presented as {\bf Theorem \ref{theo:3}}. 
The corresponding eigenvalues of the matrices are also presented as {\bf Theorem \ref{EDrelation}}.
They are polynomials in the  parameters $g$ and $h$,  
comprising of the energies of the virtual states which constitute the multi-index. 
The basic equations \eqref{basic} are `symmetric' between the original polynomial 
$P_{\mathcal{D},\,\mathcal{N}}(\eta)$ and the extra `lower degree polynomial' 
$P_{\mathcal{D}',\,\mathcal{N}}(\eta)$.  The basic equations with the roles of 
$\mathcal{D}$ and $\mathcal{D}'$ exchanged are presented as {\bf Theorem \ref{theo:6}}.
 The discrete orthogonality relations for another group of `new' orthogonal polynomials 
 are discussed in section \ref{sec:KA}. 
 They are called the Krein-Adler (K-A) polynomials based on the Hermite, Laguerre 
 and Jacobi polynomials. The above 
 {\bf Theorems \ref{theo1},\,\ref{theo:1},\,\ref{theo:2},\,\ref{theo:3},\,\ref{theo:6}} also hold for them 
  with minor modifications.
The final section is for a summary and comments.
\section{Multi-indexed Laguerre and Jacobi polynomials}
\label{sec:mulLJ}
The basic properties of the multi-indexed Laguerre and Jacobi polynomials are recapitulated here
\cite{os18,os25}.
The quantum mechanical system (Schr\"odinger equations) 
having  these polynomials as the main parts of the eigenfunctions 
\begin{equation}
  \mathcal{H}\phi_n(x)=\mathcal{E}(n)\phi_n(x)\ \ (n=0,1,\ldots),\quad
  \mathcal{H}=-\partial_x^2+U(x),
  \label{scheq}
\end{equation}
is  iso-spectrally deformed (see \eqref{defD}--\eqref{eigpro}) to
\begin{equation}
  \mathcal{H}_{\mathcal{D}}\phi_{\mathcal{D},\,n}(x)
  =\mathcal{E}(n)\phi_{\mathcal{D},\,n}(x)\ \ (n=0,1,\ldots),\quad
  \mathcal{H}_{\mathcal{D}}=-\partial_x^2+U_{\mathcal{D}}(x),
  \label{Dscheq}
\end{equation}
which gives rise to multi-indexed Laguerre and Jacobi polynomials, \eqref{XiD}--\eqref{Pn}, 
as the main parts of the eigenfunctions.
Multiple Darboux transformations \cite{darboux} are employed for the deformations 
by using two types (type I, II) of seed solutions, 
\eqref{L1seed},\eqref{L2seed},\eqref{J1seed},\eqref{J2seed}, 
called virtual state wavefunctions, 
which are obtained from the eigenfunctions through discrete symmetry transformations.

Similarly, the K-A systems are the iso-spectral deformations of the Hermite, Laguerre and Jacobi
systems  through multiple Darboux transformations.

\subsection{Original systems, eigenfunctions  and seed solutions}
\label{sec:ori}
\subsubsection{Hermite system}
\label{sec:oriher}
The quantum mechanical system corresponding to the Hermite polynomials is 
the well-known harmonic oscillator system with the quadratic potential
\begin{equation}
  U(x)\eqdef x^2-1,\quad -\infty<x<\infty.
  \label{harm}
\end{equation}
The main part of the eigenfunctions are the Hermite polynomials, $H_n(x)$,
\begin{equation}
\phi_n(x)\eqdef\phi_0(x)H_n\bigl(\eta(x)\bigr),\quad
  \mathcal{E}(n)=2n,\quad \eta(x)\eqdef x, \quad \phi_0(x)=e^{-x^2/2}.
\end{equation}
The constant ($-1$) term in the potential is added to make the groundstate energy vanishing
$\mathcal{E}(0)=0.$ The  same convention is adopted for the Laguerre and Jacobi systems
shown below.
\subsubsection{Laguerre system}
\label{sec:orilag}
The  Hamitonian with the potential 
\begin{equation}
  U(x)\eqdef x^2+\frac{g(g-1)}{x^2}-2g-1,\quad 0<x<\infty,\quad g>1/2,
  \label{rad}
\end{equation}
has the Laguerre polynomials $L^{(\alpha)}_n(\eta)$ as the main part of the
eigenfunctions:
\begin{align*}
  \phi_n(x;g)&\eqdef\phi_0(x;g)L^{(g-1/2)}_n\bigl(\eta(x)\bigr),\quad
  \mathcal{E}(n)\eqdef 4n,\quad\eta(x)\eqdef x^2,\quad
  \phi_0(x;g)\eqdef e^{-x^2/2}x^g.
\end{align*}
There are two types of seed solutions:
\begin{align} 
 \text{L1}:&\ \ \mathcal{H}\tilde{\phi}^{\I}_{\text{v}}(x)
  =\tilde{\mathcal{E}}^{\I}(\text{v})\tilde{\phi}^{\I}_{\text{v}}(x),\quad
  \tilde{\phi}^{\I}_{\text{v}}(x)\eqdef
  \tilde{\phi}^{\I}_0(x)\xi^{\I}_{\text{v}}\bigl(\eta(x)\bigr),\quad
  \tilde{\phi}^{\I}_0(x)\eqdef e^{x^2/2}x^g,\n
  &\ \ \xi^{\I}_{\text{v}}(\eta)\eqdef L^{(g-1/2)}_{\text{v}}(-\eta),
  \quad\tilde{\mathcal{E}}^{\I}(\text{v})\eqdef-4(g+\text{v}+1/2)<0,\quad
  \text{v}\in\mathbb{Z}_{>0},
   \label{L1seed}\\
    \text{L2}:&\ \ \mathcal{H}\tilde{\phi}^{\II}_{\text{v}}(x)
  =\tilde{\mathcal{E}}^{\II}(\text{v})\tilde{\phi}^{\II}_{\text{v}}(x),\quad
  \tilde{\phi}^{\II}_{\text{v}}(x)\eqdef
  \tilde{\phi}^{\II}_0(x)\xi^{\II}_{\text{v}}\bigl(\eta(x)\bigr),\quad
  \tilde{\phi}^{\II}_0(x)\eqdef e^{-x^2/2}x^{1-g},\n
  &\ \ \xi^{\II}_{\text{v}}(\eta)\eqdef L^{(1/2-g)}_{\text{v}}(\eta),
  \quad\tilde{\mathcal{E}}^{\II}(\text{v})\eqdef-4(g-\text{v}-1/2)<0,\quad
  \text{v}=1,\ldots,[g-\tfrac12]',
\label{L2seed}
\end{align}
in which $[a]'$ denotes the greatest integer less than $a$.

\subsubsection{Jacobi system}
\label{sec:orijac}
The Hamiltonian with the  potential
\begin{equation}
  U(x)\eqdef\frac{g(g-1)}{\sin^2x}+\frac{h(h-1)}{\cos^2x}-(g+h)^2,\quad
  0<x<{\pi}/{2},\quad g>1/2,\ h>1/2,
 \label{PTpot}
\end{equation}
has the Jacobi polynomials $P_n^{(\alpha,\beta)}(\eta)$ as the main part of the
eigenfunctions:
\begin{align*}
  &\phi_n\bigl(x;g,h\bigr)\eqdef\phi_0\bigl(x;g,h\bigr)
  P_n^{(g-1/2,h-1/2)}\bigl(\eta(x)\bigr),\quad
  \eta(x)\eqdef\cos2x,\\[2pt]
  &\phi_0\bigl(x;g,h\bigr)\eqdef(\sin x)^g(\cos x)^h,\quad
  \mathcal{E}(n)\eqdef 4n(n+g+h).
  \end{align*}
 Two types of seed solutions are
 \begin{align}
  &\text{J1}:\ \ \mathcal{H}\tilde{\phi}^{\I}_{\text{v}}(x)
  =\tilde{\mathcal{E}}^{\I}(\text{v})\tilde{\phi}^{\I}_{\text{v}}(x),\quad
  \tilde{\phi}^{\I}_{\text{v}}(x)\eqdef
  \tilde{\phi}^{\I}_0(x)\xi^{\I}_{\text{v}}\bigl(\eta(x)\bigr),\quad
  \tilde{\phi}^{\I}_0(x)\eqdef(\sin x)^g(\cos x)^{1-h},\n
  &\ \xi^{\I}_{\text{v}}(\eta)\eqdef
  P^{(g-1/2,1/2-h)}_{\text{v}}(\eta),
  \ \tilde{\mathcal{E}}^{\I}(\text{v})\eqdef
  -4(g+\text{v}+1/2)(h-\text{v}-1/2)<0,
  \  \text{v}=1,\ldots,[h-\tfrac12]',
  \label{J1seed}\\
 & \text{J2}:\ \ \mathcal{H}\tilde{\phi}^{\II}_{\text{v}}(x)
  =\tilde{\mathcal{E}}^{\II}(\text{v})\tilde{\phi}^{\II}_{\text{v}}(x),\quad
  \tilde{\phi}^{\II}_{\text{v}}(x)\eqdef
  \tilde{\phi}^{\II}_0(x)\xi^{\II}_{\text{v}}\bigl(\eta(x)\bigr),\quad
  \tilde{\phi}^{\II}_0(x)\eqdef(\sin x)^{1-g}(\cos x)^h,\n
  &\ \xi^{\II}_{\text{v}}(\eta)\eqdef
  P^{(\frac12-g,h-1/2)}_{\text{v}}(\eta),
  \ \tilde{\mathcal{E}}^{\II}(\text{v})\eqdef
  -4(g-\text{v}-1/2)(h+\text{v}+1/2)<0,
  \  \text{v}=1,\ldots,[g-\tfrac12]'.
  \label{J2seed}
\end{align}

\subsection{Wronskians for the multi-indexed Laguerre and Jacobi
polynomials}
\label{sec:wron}
The deformed systems are obtained by multiple Darboux transformations in terms of two types of 
seed solutions specified by the multi-index $\mathcal{D}$ consisting of  positive
integers which are the {\bf degrees} of the virtual state wavefunctions,

\begin{align}
&\mathcal{D}=\mathcal{D}^{\I}\cup\mathcal{D}^{\II},\quad
\mathcal{D}^{\I}=\{d_1^{\I},\ldots,d_M^{\I}\},\qquad
\mathcal{D}^{\II}=\{d_1^{\II},\ldots,d_N^{\II}\}, 
\label{defD}\\
&\phantom{\mathcal{D}=\mathcal{D}^{\I}\cup\mathcal{D}^{\II},\quad}
0<d_1^{\I}<\cdots<d_M^{\I},
\qquad\
0<d_1^{\II}<\cdots<d_N^{\II},\n
  U_{\mathcal{D}}(x)&\eqdef U(x)-2\partial_x^2
  \log\bigl|\text{W}[\tilde{\phi}_{d_1^{\I}}^{\I},\ldots,\tilde{\phi}_{d_M^{\I}}^{\I},
  \tilde{\phi}_{d_1^{\II}}^{\II},\ldots,\tilde{\phi}_{d_N^{\II}}^{\II}](x)\bigr|,\\
  \phi_{\mathcal{D},\,n}(x)&\eqdef
  \frac{\text{W}[\tilde{\phi}_{d_1^{\I}}^{\I},\ldots,\tilde{\phi}_{d_M^{\I}}^{\I},
  \tilde{\phi}_{d_1^{\II}}^{\II},\ldots,\tilde{\phi}_{d_N^{\II}}^{\II},\phi_n](x)}
  {\text{W}[\tilde{\phi}_{d_1^{\I}}^{\I},\ldots,\tilde{\phi}_{d_M^{\I}}^{\I},
  \tilde{\phi}_{d_1^{\II}}^{\II},\ldots,\tilde{\phi}_{d_N^{\II}}^{\II}](x)}
  \ \ (n=0,1,\ldots)\\
  &\eqdef\cF^{M+N}\psi_{\mathcal{D}}(x)
  P_{\mathcal{D},n}\bigl(\eta(x)\bigr),\quad
  \psi_{\mathcal{D}}(x)\eqdef\frac{\hat{\phi}_0(x)}
  {\Xi_{\mathcal{D}}\bigl(\eta(x)\bigr)},\\
\mathcal{H}_{\mathcal D} \phi_{\mathcal{D},\,n}(x)&
=\mathcal{E}(n) \phi_{\mathcal{D},\,n}(x),\qquad 
\mathcal{H}_{\mathcal D}=-\frac{d^2}{dx^2}+U_{\mathcal D}(x),
\quad n=0,1,\ldots,
\label{eigpro}
\end{align}
in which $\hat{\phi}_0(x)$ and $\cF$ are defined by
\begin{equation*}
  \hat{\phi}_0(x)\eqdef\left\{
  \begin{array}{ll}
  \phi_0(x;g+M-N)&:\text{L}\\[2pt]
  \phi_0\bigl(x;g+M-N,h-M+N\bigr)&:\text{J}
  \end{array}\right.,\quad
  \cF\eqdef\left\{
  \begin{array}{cl}
  2&:\text{L}\\
  -4&:\text{J}
  \end{array}\right..
\end{equation*}
The denominator polynomial $\Xi_{\mathcal{D}}(\eta)$ and the multi-indexed
orthogonal polynomial $P_{\mathcal{D},\,n}(\eta)$ are also expressed by Wronskians,
 \begin{align}
  \Xi_{\mathcal{D}}(\eta)
  &\eqdef\text{W}[\mu_{d_1^{\I}},\ldots,\mu_{d_N^{\II}}](\eta)\times\left\{
  \begin{array}{ll}
  \eta^{(M+g-\frac12)N}e^{-M\eta}&:\text{L}\\
  \bigl(\frac{1-\eta}{2}\bigr)^{(M+g-\frac12)N}
  \bigl(\frac{1+\eta}{2}\bigr)^{(M+h-\frac12)N}&:\text{J}
  \end{array}\right.,
  \label{XiD}\\
  P_{\mathcal{D},n}(\eta)
  &\eqdef\text{W}[\mu_{d_1^{\I}},\ldots,\mu_{d_N^{\II}},P_n](\eta)\times\left\{
  \begin{array}{ll}
  \eta^{(M+g+\frac12)N}e^{-M\eta}&:\text{L}\\
  \bigl(\frac{1-\eta}{2}\bigr)^{(M+g+\frac12)N}
  \bigl(\frac{1+\eta}{2}\bigr)^{(N+h+\frac12)N}&:\text{J}
  \end{array}\right.,
  \label{PDn}\\
  \mu_{\text{v}}(\eta)&\eqdef\left\{
  \begin{array}{ll}
  e^{\eta}\times L^{(g-\frac12)}_{\text{v}}(-\eta)
  &:\text{L, $\text{v}$ type $\I$}\\
  \eta^{\frac12-g}\times L^{(\frac12-g)}_{\text{v}}(\eta)
  &:\text{L, $\text{v}$ type $\II$}\\
  \bigl(\frac{1+\eta}{2}\bigr)^{\frac12-h}\times
  P^{(g-\frac12,\frac12-h)}_{\text{v}}(\eta)
  &:\text{J, $\text{v}$ type $\I$}\\[3pt]
  \bigl(\frac{1-\eta}{2}\bigr)^{\frac12-g}\times
  P^{(\frac12-g,h-\frac12)}_{\text{v}}(\eta)
  &:\text{J, $\text{v}$ type $\II$}
  \end{array}\right.,
  \label{muv}\\
  P_n(\eta)
  &\eqdef\left\{
  \begin{array}{ll}
  L^{(g-\frac12)}_n(\eta)&:\text{L}\\[2pt]
  P^{(g-\frac12,h-\frac12)}_n(\eta)&:\text{J}
  \end{array}\right..
  \label{Pn}
\end{align}
The special cases of either $M=0$ (type II only) or $N=0$ (type I only) are meaningful.
For the non-singularity of the deformed potential $U_{\mathcal D}(x)$, 
the virtual state energies $\tilde{\mathcal E}^{\ai,\ait}(\text{v})$ must be negative and
this imposes certain bounds on the parameters $g$ and $h$:
\begin{align}
  \text{L}:&\quad
  g>\text{max}\{N+3/2,d_j^\ait+1/2\},
  \label{Lbound}\\
  \text{J}:&\quad
  g>\text{max}\{N+2,d_j^\ait+1/2\},\quad
  h>\text{max}\{M+2,d_j^\ai+1/2\}.
  \label{Jbound}
\end{align}

Although we have restricted $d_j^\I, d_j^\II$ to be positive, including the case of $d_j^\I, d_j^\II=0$ 
goes almost parallel.  
However, inclusion of the latter does not give new polynomials because 
of the shape-invariance \cite{os25}. 

For generic values of the  parameters $g$ and $h$, the multi-indexed polynomial 
$P_{\mathcal{D},n}(\eta)$ is of degree $\ell_{\mathcal D}+n$
and the denominator polynomial $\Xi_{\mathcal D}(\eta)$ is of degree $\ell_{\mathcal D}$ in
$\eta$, in which $\ell_{\mathcal D}$ is given by
\begin{equation}
\ell_{\mathcal D}\eqdef \ell_{{\mathcal D}^\ai}+\ell_{{\mathcal D}^\ait}+MN\ge1,\quad
  \ell_{{\mathcal D}^\ai}\eqdef\sum_{j=1}^Md_j^\ai -\frac12 M(M-1),\quad
   \ell_{{\mathcal D}^\ait}\eqdef\sum_{j=1}^Nd_j^\ait
 -\frac12N(N-1).
 \label{ellD}
\end{equation}
It may happen that for some exceptional values of the parameters, the actual degrees  
are smaller than $\ell_{\mathcal D}$  given above.

The second order differential operator
$\widetilde{\mathcal{H}}_{\mathcal{D}}$ governing the
multi-indexed polynomial $P_{\mathcal{D},\,n}\,(\eta)$ is obtained from $\mathcal{H}_{\mathcal D}$ 
by a similarity transformation in terms of $\psi_{\mathcal D}(x)$:
\begin{align}
  \widetilde{\mathcal{H}}_{\mathcal{D}}  &\eqdef\psi_{\mathcal{D}}(x)^{-1}\circ
  \mathcal{H}_{\mathcal{D}}\circ
  \psi_{\mathcal{D}}(x)=-4\left(c_2(\eta)\frac{d^2}{d\eta^2}
  +f_1(\eta)\frac{d}{d\eta}+f_2(\eta)\right), 
   \label{ThamD}\\
 &\qquad f_1(\eta)\eqdef 
c_1(\eta,\bm{\lambda}^{[M,N]})-2c_2(\eta)
  \frac{\partial_{\eta}\Xi_{\mathcal{D}}(\eta)}
  {\Xi_{\mathcal{D}}(\eta)},\\
 &\qquad f_2(\eta)\eqdef
 c_2(\eta)
  \frac{\partial^2_{\eta}\Xi_{\mathcal{D}}(\eta)}
  {\Xi_{\mathcal{D}}(\eta)}
  -c_1(\eta,\bm{\lambda}^{[M,N]}-\bm{\delta})
  \frac{\partial_{\eta}\Xi_{\mathcal{D}}(\eta)}
  {\Xi_{\mathcal{D}}(\eta)},
  \label{f2}\\
  &\hspace{-14mm}\widetilde{\mathcal{H}}_{\mathcal{D}}
  P_{\mathcal{D},\,n}(\eta)=\mathcal{E}(n)
  P_{\mathcal{D},\,n}(\eta),
  \label{fuchs}
\end{align}
in which $\bm{\lambda}$ stands for the  parameters, $\bm{\lambda}=g$ for L and 
$\bm{\lambda}=(g,h)$ for J and 
\begin{equation}
  c_1(\eta,\bm{\lambda})\eqdef\left\{
  \begin{array}{ll}
  g+\tfrac12-\eta&\!:\text{L}\\
  h-g-(g+h+1)\eta&\!:\text{J}
  \end{array}\right.\!\!,
  \ \ c_2(\eta)\eqdef\left\{
  \begin{array}{ll}
  \eta&\!:\text{L}\\
  1-\eta^2&\!:\text{J}
  \end{array}\right..
  \label{cF,c1,c2}
\end{equation}
The shifted parameters are
\begin{equation}
  \bm{\lambda}^{[M,N]}=g+M-N\ \ \text{for L},\quad
  \bm{\lambda}^{[M,N]}=(g+M-N,h-M+N)\ \ \text{for J},
  \label{parashift}
\end{equation}
and the shift $\bm{\delta}$ is $\bm{\delta}=1$
for L and $\bm{\delta}=(1,1)$ for J.

It should be stressed that  the second order differential operator 
$\Xi_{\mathcal{D}} \widetilde{\mathcal{H}}_{\mathcal{D}}$ maps any polynomial (in $\eta$)
to  a polynomial.

\section{Discrete orthogonality relations for  multi-indexed Laguerre and Jacobi polynomials}
\label{sec:disortmulLJ}

In this section we answer the three questions {\bf Q1}--{\bf Q3} affirmatively.
The discrete orthogonality relations hold definitely for the multi-indexed 
Laguerre and Jacobi polynomials ({\bf Q1}),  the identification of the extra `lower degree polynomials'
is quite natural ({\bf Q2}) and the  Christoffel numbers
have the same explicit forms as those for the original Laguerre and Jacobi polynomials ({\bf Q3}).

Before presenting the main results, let us review the discrete orthogonality relations for the classical
orthogonal polynomials, the Hermite (H), Laguerre (L) and Jacobi (J) polynomials, 
which were obtained in \cite{zero}:
\begin{align} 
\text{H:} & \quad \sum_{l=1}^\mathcal{N}\frac{H_n(x_l)H_m(x_l)}
{\bigl(H'_{\mathcal N}(x_l)\bigr)^2}=0,\quad 
  m\neq n, \quad m,n=0,1,\ldots, \mathcal{N}-1,
  \label{disH}\\
\text{L:} & \quad \sum_{l=1}^\mathcal{N}\frac{L_n^{(\alpha)}(\eta_l)L_m^{(\alpha)}(\eta_l)}
{\eta_l\bigl({L^{(\alpha)}_{\mathcal N}}'(\eta_l)\bigr)^2}=0,\quad 
  m\neq n, \quad m,n=0,1,\ldots, \mathcal{N}-1,\\
\text{J:} & \quad \sum_{l=1}^\mathcal{N}
\frac{P_n^{(\alpha,\beta)}(\eta_l)P_m^{(\alpha,\beta)}(\eta_l)}
{(1-\eta_l^2)\bigl({P^{(\alpha,\beta)}_{\mathcal N}}'(\eta_l)\bigr)^2}=0,\quad 
  m\neq n, \quad m,n=0,1,\ldots, \mathcal{N}-1.  
  \label{disJ}
\end{align}
These formulas, including the other members of the classical orthogonal polynomials in the Askey scheme 
\cite{koeswart}, can be expressed in a unified fashion.
Let us express the classical orthogonal polynomials by using the sinusoidal coordinate $\eta(x)$ as
 $\{P_n\bigl(\eta(x)\bigr)\}$. 
 As shown above $\eta(x)=x$ for H, $\eta(x)=x^2$ for L and $\eta(x)=\cos2x$ for J.
 For the definitions of $\eta(x)$ for the other Askey scheme polynomials, see \cite{zero}.
 For a given positive  integer $\mathcal{N}$, let us denote the zeros of $P_{\mathcal N}\bigl(\eta(x)\bigr)$ as
 \begin{equation*}
\eta_1\equiv \eta(x_1),\ \eta_2\equiv \eta(x_2),\ \ldots,\ 
\eta_{\mathcal N}\equiv \eta(x_{\mathcal N}).
\end{equation*}
The universal discrete orthogonality relations read
\begin{equation}
\sum_{l=1}^\mathcal{N}\frac{P_n(\eta_l)P_m(\eta_l)}{\dot{\eta}(x_l)^2\bigl({P_{\mathcal N}}'(\eta_l)\bigr)^2}=0,\quad 
  m\neq n, \quad m,n=0,1,\ldots, \mathcal{N}-1,
  \label{univort}
\end{equation}
in which 
\begin{equation}
\dot{\eta}(x)\eqdef \frac{d\eta(x)}{dx}.
\end{equation}
In other words, the Christoffel numbers have a universal expression for the above  
classical orthogonal polynomials:
\begin{equation}
\lambda_j\propto \frac1{\dot{\eta}(x_j)^2\bigl({P_{\mathcal N}}'(\eta_j)\bigr)^2}>0, 
\quad j=1,\ldots, \mathcal{N}.
\label{univC}
\end{equation}
Note that $\dot{\eta}(x_j)^2$ is always expressed by $\eta_j\equiv \eta(x_j)$,
{\em e.g.\/} $\dot{\eta}(x_j)^2=4\eta_j$ for L and $\dot{\eta}(x_j)^2=4(1-\eta_j^2)$ for J.
The best known example is the Chebyshev polynomial of the first kind, 
$P_{\mathcal N}\bigl(\eta(x)\bigr)\eqdef \cos \mathcal{N}x$, $\eta(x)=\cos x$. 
In this case the Christoffel numbers are
constant, as $\dot{\eta}(x)^2=\sin^2x$, $\sin^2xP_{\mathcal N}'(\eta)^2
=\mathcal{N}^2\sin^2\mathcal{N}x=\mathcal{N}^2$, 
at $x_l$ such that $\cos\mathcal{N}x_l=0$.
As we will see in the subsequent sections, this universal formula applies to the multi-indexed 
Laguerre and  Jacobi polynomials and the K-A orthogonal polynomials, too 
({\it cf.}  {\bf Corollary \ref{chris4}}).
The above formula \eqref{univC}  is different from the general formula reported 
in the textbooks for orthogonal polynomials,
see for example (5.6.7) in \cite{askey}.

Several interesting examples of other types of orthogonality relations for the  Laguerre and Jacobi 
polynomials and for the Bessel functions
were reported in \cite{ahmed},  see  (4.15) for L , (5.15), (5.25) for J and (6.27) for B.
\subsection{Extra `lower degree polynomials' }
\label{sec:lowdeg}

Let us now investigate whether the discrete orthogonality relations hold or not 
for the multi-indexed Laguerre and Jacobi polynomials.
We fix a positive integer $\mathcal{N}$ and 
try to formulate the orthogonality relations like \eqref{univort}
in terms of the zeros of $P_{\mathcal{D},\,\mathcal{N}}\bigl(\eta(x)\bigr)$;
\begin{equation}
P_{\mathcal{D},\,\mathcal{N}}\bigl(\eta(x)\bigr)=0,\quad
\eta_1\equiv \eta(x_1),\ \eta_2\equiv \eta(x_2),\ \ldots,
\ \eta_{\widetilde{\mathcal{N}}}\equiv \eta(x_{\widetilde{\mathcal{N}}}),\quad
\widetilde{\mathcal{N}}\eqdef\ell_{\mathcal D}+\mathcal{N}.
\label{lDNzeros}
\end{equation}
As far as we know there are four references discussing the zeros of the exceptional (one-indexed) Laguerre and 
Jacobi polynomials \cite{hszero,gmmzero,bkzero,bzero}.
The number of polynomials participating  in the discrete orthogonality relations 
is the same as the number of the zeros.
This means we have to {\em identify the extra $\ell_{\mathcal D}$ 
`lower degree polynomials' \/} on top of the 
ordinary ones,
\begin{equation*}
\{P_{\mathcal{D},\,m}(\eta)\} \quad m=0,1,\ldots, \mathcal{N}-1.
\end{equation*}

A hint to solve this problem comes from the Wronskian expression \eqref{PDn} for the multi-indexed 
Laguerre and Jacobi polynomials consisting of seed polynomials of degree $d^\ai_j$ ($j=1,\ldots,M$)
and $d^\ait_k$ ($k=1,\ldots,N$) and the original polynomial of degree $\mathcal{N}$.
The above ordinary lower degree polynomials are obtained by reducing $\mathcal{N}\to m$.
A good guess and a correct one is to reduce the degrees of one of the seed polynomials,
$d^\ai_j$ or $d^\ait_k$. 
As this formulation is the same for type I and II, we explain it  for type I by omitting the superscript I.
Let $F$ be the set of all degrees up to $d_M=\text{max}(\mathcal{D})$, the maximum degree, 
and we define $E$ as the complement set
of $\mathcal{D}$ in $F$, specifying the empty (not used) degrees:
\begin{equation}
\mathcal{D}=\{d_1,\ldots,d_M\}, \quad
F\eqdef\{0,1,\ldots,d_M\},\quad
E\eqdef F\backslash\mathcal{D}.
\label{defE}
\end{equation}
For each $d_j\in\mathcal{D}$ let us define the set of lower degrees as
\begin{equation}
E_j\eqdef \{\epsilon\in E| \epsilon<d_j\},\quad j=1,\ldots,M.
\label{edef}
\end{equation}
By replacing $d_j$ with one of $\epsilon\in E_j$ in the Wronskian expression \eqref{PDn} 
for $P_{\mathcal{D},\,\mathcal{N}}(\eta)$, we obtain a
lower degree ($\ell_{\mathcal D}+\mathcal{N}-d_j^\ai+\epsilon_k^\ai$) polynomial, 
which is specified by a multi-index
\begin{align}
\mathcal{D}'_{\ai\,j,k} \eqdef  {\mathcal D}^{\prime\ai}_{j,k}\cup \mathcal{D}^\ait,\quad
{\mathcal D}^{\prime\ai}_{j,k} \eqdef 
\left(\mathcal{D}^\ai \backslash \{d_j^\ai\}\right)\cup \{\epsilon_k^\ai\}
=\{d_1^\ai,\ldots,d_{j-1}^\ai,\epsilon_k^\ai,d_{j+1}^\ai,\ldots,d_M^\ai\},\n
~ d_j^\ai\in {\mathcal D}^{\ai},
~ \epsilon_k^\ai\in E_j^\ai,
\label{1D}
\end{align}
in which we use $\epsilon_k^\ai$ to show that $\epsilon$ belongs to $E_j^\ai$.
Let us refer to these extra lower degree polynomials as type I lower degree polynomials.
The total number of such lower degree polynomials is 
\begin{equation*}
\sum_{j=1}^M\left(d_j^\ai-(j-1)\right)=\sum_{j=1}^Md_j^\ai-\tfrac12M(M-1)=\ell_{D^\ai}.
\end{equation*}

Similarly we construct a  type II lower degree 
($\ell_{\mathcal D}+\mathcal{N}-d_j^\ait+\epsilon_k^\ait$)  
polynomial specified by a  multi-index
\begin{align}
\mathcal{D}'_{\ait\,j,k}\eqdef 
\mathcal{D}^\ai\cup  {\mathcal D}^{\prime\ait}_{j,k},\quad   
{\mathcal D}^{\prime\ait}_{j,k}\eqdef 
\left(\mathcal{D}^\ait \backslash \{d_j^\ait\}\right)\cup \{\epsilon_k^\ait\}
=\{d_1^\ait,\ldots,d_{j-1}^\ait,\epsilon_k^\ait,d_{j+1}^\ait,\ldots,d_N^\ait\},\n
~ d_j^\ait\in {\mathcal D}^\ait,
~ \epsilon_k^\ait\in E_j^\ait. 
\label{2D}
\end{align}
There are $\ell_{\mathcal{D}^\ait}$ such polynomials.
The member of the third group of $MN$ (see \eqref{ellD}) lower degree polynomials is obtained by
removing $d_j^\ai\in\mathcal{D}^\ai$ and $d_k^\ait\in\mathcal{D}^\ait$ simultaneously. 
It has a degree $\ell_{\mathcal D}+\mathcal{N}-d_j^\ai-d_k^\ait-1$ 
corresponding to the multi-index
\begin{align}
\mathcal{D}'_{\at\,j,k}&\eqdef \{d_1^\ai,\ldots,d_{j-1}^\ai,d_{j+1}^\ai,\ldots,d_M^\ai\}\cup
\{d_1^\ait,\ldots,d_{k-1}^\ait,d_{k+1}^\ait,\ldots,d_N^\ait\}.
\label{3D}
\end{align}
The removed indices are $d_j^\ai$ and $d_k^\ait$.

The set of all multi-indices corresponding to the Extra `lower degree Polynomials' for given
$\mathcal{D}$ is denoted by $\mathcal{EP}$, $\#\mathcal{EP}=\ell_{\mathcal D}$.
\subsection{Discrete orthogonality relations}
\label{sec:discort}
Let us denote the members of the extra lower degree polynomials symbolically as 
$P_{\mathcal{D}',\,\mathcal{N}}(\eta)$ and $P_{\mathcal{D}'',\,\mathcal{N}}(\eta)$, etc.
We have the following
\begin{theo}
\label{theo1}
For the multi-indexed Laguerre and Jacobi polynomials, 
the discrete orthogonality relations take the following forms;
\begin{align} 
\sum_{l=1}^{\widetilde{\mathcal{N}}}
\frac{P_{\mathcal{D},\,n}(\eta_l)P_{\mathcal{D},\,m}(\eta_l)}
{\dot{\eta}(x_l)^2\bigl(P_{\mathcal{D},\,\mathcal{N}}'(\eta_l)\bigr)^2}&=0,\quad 
  m\neq n, \quad m,n=0,1,\ldots, \mathcal{N}-1,
 \label{disor1}\\
\sum_{l=1}^{\widetilde{\mathcal{N}}}
\frac{P_{\mathcal{D},\,m}(\eta_l)P_{\mathcal{D}',\,\mathcal{N}}(\eta_l)}
{\dot{\eta}(x_l)^2\bigl(P_{\mathcal{D},\,\mathcal{N}}'(\eta_l)\bigr)^2}&=0,\quad 
 m=0,1,\ldots, \mathcal{N}-1,\quad \mathcal{D}'  \subseteq\mathcal{EP},
 \label{disor2}\\
\sum_{l=1}^{\widetilde{\mathcal{N}}}
\frac{P_{\mathcal{D}',\,\mathcal{N}}(\eta_l)P_{\mathcal{D}'',\,\mathcal{N}}(\eta_l)}
{\dot{\eta}(x_l)^2\bigl(P_{\mathcal{D},\,\mathcal{N}}'(\eta_l)\bigr)^2}&=0,\quad 
\mathcal{D}'\neq\mathcal{D}'',\quad   \mathcal{D}',\mathcal{D}'' \subseteq \mathcal{EP}.
 \label{disor3}
\end{align}
\end{theo}
In the present case, the universal Christoffel numbers \eqref{univC} are no longer positive. 
Similar discrete orthogonality relations are expected to hold for the multi-indexed 
Racah, $q$-Racah, Wilson and Askey-Wilson polynomials \cite{dismul}. 
Those generated by Krein-Adler transformations
\cite{krein,adler} will be discussed in section \ref{sec:KA}.

\begin{coro}
The polynomials participating in the above discrete orthogonality relations are linearly independent.
The above extra `lower degree polynomials' and the multi-indexed polynomials 
$\{P_{\mathcal{D},\,n}(\eta)\}$, $n=0,\ldots,\mathcal{N}$ form the basis of the space of polynomials 
up to degree $\widetilde{\mathcal N}=\ell_{\mathcal D}+\mathcal{N}$.
\end{coro}

\bigskip
We will present a proof in the subsequent section.
The explicit forms of the discrete orthogonality relations for the classical orthogonal polynomials 
\eqref{disH}--\eqref{univort} were obtained as a consequence of the orthogonality among
the eigenvectors of a certain hermitian (real symmetric) matrix $\mathcal{M}$.
The matrix was constructed within the quantum mechanical formulation 
of classical orthogonal polynomials through perturbation around the zeros of truncated wavefunctions
at degree $\mathcal{N}$  \cite{zero}. 
In the subsequent section we apply the same method to the multi-indexed
Laguerre and Jacobi polynomials to gain further insights from different angles, {\em e.g.\/} 
the rationale for the orthogonality, the eigenvalues corresponding to the extra `lower degree polynomials',
etc.

\section{Perturbation around the zeros}
\label{sec:zerosum}

The starting point is the  time-dependent one-dimensional Schr\"odinger equation
\begin{equation}
i\frac{\partial \Psi(x,t)}{\partial t}=\mathcal{H}_{\mathcal D} \Psi(x,t),
\label{scheq2}
\end{equation}
for the quantum mechanics of the multi-indexed Laguerre and Jacobi polynomials \eqref{Dscheq}.
The Hamiltonian $\mathcal{H}_{\mathcal D}$ \eqref{eigpro} is  time independent.
In terms of the complete set of solutions $\{\mathcal{E}(n),\phi_{\mathcal{D},\,n}(x)\}$ of 
the eigenvalue problem \eqref{eigpro}, the general solution of the time-dependent
Schr\"odinger equation \eqref{scheq2} is given by
\begin{equation}
\Psi(x,t)=\sum_{n=0}^\infty c_n e^{-i\mathcal{E}(n)t}\phi_{\mathcal{D},\,n}(x),
\label{gensol}
\end{equation}
in which $\{c_n\}$ are the constants of integration.
As shown in section \ref{sec:mulLJ} the eigenfunctions have a factorised form
$\phi_{\mathcal{D},\,n}(x)
  \propto\psi_{\mathcal{D}}(x)
 P_{\mathcal{D},\,n}(\eta(x))$,
and  the multi-indexed polynomials $\{P_{\mathcal{D},\,n}(\eta(x))\}$ 
satisfy second order differential equations \eqref{ThamD}--\eqref{fuchs},
\begin{equation*}
\widetilde{\mathcal H}_{\mathcal{D}}P_{\mathcal{D},\,n}\left(\eta(x)\right)
=\mathcal{E}(n)P_{\mathcal{D},\,n}\left(\eta(x)\right),\quad
n=0,1,\ldots .
\end{equation*}

\subsection{Polynomial solutions}
\label{sec:poly}

Let us fix a positive integer $\mathcal{N}$ 
and restrict the general solution \eqref{gensol}
to those having {\em degrees up to \/} 
$\widetilde{\mathcal{N}}\eqdef \ell_{\mathcal{D}}+\mathcal{N}$:
\begin{align} 
 \Psi_{\widetilde{\mathcal{N}}}(x,t)&=\sum_{n=0}^{\mathcal{N}} 
 c_n e^{-i\mathcal{E}(n)t}\phi_{\mathcal{D},\,n}(x),\\
 &=e^{-i\mathcal{E}(\mathcal{N})t}\psi_{\mathcal{D}}(x)\psi_{\widetilde{\mathcal N}}(x,t).
 \label{polysol}
\end{align}
Here the function
\begin{equation}
\psi_{\widetilde{\mathcal{N}}}(x,t)\eqdef \sum_{n=0}^{\mathcal N}c_n' 
e^{i(\mathcal{E}(\mathcal{N})-\mathcal{E}(n))t}P_{\mathcal{D},\,n}\left(\eta(x)\right),
\quad c_n'\eqdef c_n{\cF}^{M+N},
\label{psiNdef}
\end{equation}
is a polynomial of degree $\widetilde{\mathcal{N}}$ in $\eta(x)$.
We choose the coefficient $c_{\mathcal N}'$ of the highest degree polynomial  
$P_{\mathcal{D},\,\mathcal N}\left(\eta(x)\right)$
to make it monic:
\begin{equation*}
c_{\mathcal N}'P_{\mathcal{D},\,\mathcal N}\bigl(\eta(x)\bigr)
=\prod_{n=1}^{\widetilde{\mathcal N}}\bigl(\eta(x)-\eta(x_n)\bigr),
\end{equation*}
in which $\{\eta(x_n)\}$, $n=1,\ldots,\widetilde{\mathcal N}$,  
are the {\em zeros} of $P_{\mathcal{D},\,\mathcal N}\left(\eta(x)\right)$ \eqref{lDNzeros}.
Since all the coefficients of $P_{\mathcal{D},\,\mathcal N}\bigl(\eta)$ are real, 
the zeros are real or complex conjugate pairs:
\begin{equation*}
\eta_n\eqdef\eta(x_n),\quad 
\{\eta_1,\eta_2,\ldots,\eta_{\widetilde{\mathcal N}}\}
=\{\eta_1^*,\eta_2^*\ldots,\eta_{\widetilde{\mathcal N}}^*\}
\quad \text{as a set}.
\end{equation*}
Among them, there are  $\mathcal{N}$ {\em real} zeros in the orthogonality domain $0<\eta<\infty$ for L
and $-1<\eta<1$ for J. Let us call them the {\em ordinary zeros} and the rest $\ell_{\mathcal D}$ zeros
as the {\em extra zeros}.
Among the extra zeros, there are real zeros and complex conjugate pairs, the latter will be denoted as
$\{\eta_n,\eta_{\bar{n}}\}$, $(\eta_n)^*=\eta_{\bar{n}}$.
We also use this notation for all real zeros {\em i.e.\/} $\bar{n}\eqdef n$ for a real zero $\eta_n$.

Obviously the polynomial $\psi_{\widetilde{\mathcal N}}(x,t)$ \eqref{psiNdef} 
satisfies the time evolution equation
\begin{align}
\frac{\partial \psi_{\widetilde{\mathcal N}}(x,t)}{\partial t}
=-i\widetilde{\mathcal H}_{\mathcal{D},\,\mathcal N}\psi_{\widetilde{\mathcal N}}(x,t),\qquad
\widetilde{\mathcal H}_{\mathcal{D},\,\mathcal N}\eqdef 
\widetilde{\mathcal H}_{\mathcal D}-\mathcal{E}({\mathcal N}).
\label{tevo}
\end{align}
As $\psi_{\widetilde{\mathcal N}}(x,t)$ is a degree $\widetilde{\mathcal{N}}$ polynomial in $\eta(x)$, 
its zeros are functions of $t$.
Let us assume that  $\psi_{\widetilde{\mathcal N}}(x,t)$ is a $t$-dependent deformation
of the highest degree monic polynomial 
$c_{\mathcal N}'P_{\mathcal{D},\, \mathcal N}\left(\eta(x)\right)$
of the following form
\begin{equation}
\psi_{\widetilde{\mathcal N}}(x,t)=
\prod_{n=1}^{\widetilde{\mathcal N}}\bigl(\eta(x)-\eta(x_n(t))\bigr),
\label{xntdef}
\end{equation}
in which $\{x_n(t)\}$ are certain  differentiable functions of $t$, describing the zeros of 
$\psi_{\widetilde{\mathcal N}}(x,t)$ at time $t$.

Among the  $t$-dependent deformations \eqref{xntdef} of the multi-indexed polynomial
$P_{\mathcal{D},\,\mathcal N}\left(\eta(x)\right)$, let us focus on those describing 
{\em infinitesimal
oscillations around the zeros\/} of $P_{\mathcal{D},\,\mathcal N}\left(\eta(x)\right)$:
\begin{equation}
x_n(t)=x_n+\epsilon\gamma_n(t),\quad 0<\epsilon\ll1,\quad 
n=1,\ldots, \widetilde{\mathcal N}.
\label{ansatz}
\end{equation}
This corresponds to choosing infinitesimal $\{c_n\}$ 
so that the deformation can be considered as perturbations around the 
zeros of $P_{\mathcal{D},\,\mathcal N}(\eta(x))$.
With the above ansatz 
$\psi_{\widetilde{\mathcal N}}(x,t)$ reads
\begin{equation*} 
\psi_{\widetilde{\mathcal N}}(x,t)=\prod_{n=1}^{\widetilde{\mathcal N}}\bigl(\eta(x)-\eta_n\bigr)
             -\epsilon\sum_{n=1}^{\widetilde{\mathcal N}}\gamma_n(t)\dot{\eta}(x_n)
    \prod_{j\neq n}^{\widetilde{\mathcal N}}\bigl(\eta(x)-\eta_j\bigr) +O(\epsilon^2),
\end{equation*}
and  each side of the time evolution equation \eqref{tevo} reduces to:
\begin{align*}
\text{l.h.s.:}&\quad-\epsilon\sum_{n=1}^{\widetilde{\mathcal N}}\frac{d\gamma_n(t)}{dt}\dot{\eta}(x_n)
           \prod_{j\neq n}^{\widetilde{\mathcal N}}\bigl(\eta(x)-\eta_j\bigr) +O(\epsilon^2),\\
\text{r.h.s.:}&\quad i\epsilon\sum_{m=1}^{\widetilde{\mathcal N}}\gamma_m(t)\dot{\eta}(x_m)
         \widetilde{\mathcal H}_{\mathcal{D},\, \mathcal N}\prod_{j\neq m}^{\widetilde{\mathcal N}}
             \bigl(\eta(x)-\eta_j\bigr) +O(\epsilon^2),
\end{align*}
since the leading polynomial $P_{\mathcal{D},\,\mathcal N}\left(\eta(x)\right)\propto 
\prod_{n=1}^{\widetilde{\mathcal N}}\left(\eta(x)-\eta_n\right)$ is annihilated by 
$\widetilde{\mathcal H}_{\mathcal{D},\,\mathcal N}$.
The l.h.s. is a degree $\widetilde{\mathcal N}-1$ polynomial in $\eta(x)$.
The r.h.s. is a rational function of $\eta(x)$ due to the denominator polynomial $\Xi_{\mathcal D}(\eta)$
in $\widetilde{\mathcal H}_{\mathcal{D},\,\mathcal{N}}$.

Let us investigate this equation at the $\widetilde{\mathcal N}$ zeros  $\{\eta(x_n)\}$ of 
$P_{\mathcal{D},\,\mathcal N}\left(\eta(x)\right)$.
This leads to $\widetilde{\mathcal N}$ linear ODE's for the unknown functions $\{\gamma_n(t)\}$
at the leading order of $\epsilon$:
\begin{align*} 
&\frac{d\gamma_n(t)}{dt}\dot{\eta}(x_n) 
\prod_{j\neq n}^{\widetilde{\mathcal N}}\bigl(\eta_n-\eta_j\bigr)\\
\quad &=-i\sum_{m=1}^{\widetilde{\mathcal N}}\gamma_m(t)\dot{\eta}(x_m)
\left.\left(\widetilde{\mathcal H}_{\mathcal{D},\, \mathcal N}
\prod_{j\neq m}^{\mathcal N}\bigl(\eta(x)-\eta_j\bigr)
\right)\right |_{x=x_n},\quad n=1,\ldots,\widetilde{\mathcal N},
\end{align*}
which can be rewritten in a matrix form:
\begin{align}   
\frac{d\gamma_n(t)}{dt}&= i\sum_{m=1}^{\widetilde{\mathcal N}}\mathcal{M}_{n\,m}\gamma_m(t), 
\quad  n=1,\ldots,\widetilde{\mathcal N}, 
\label{mateq}\\
\mathcal{M}_{n\,m}&\eqdef -\frac{\dot{\eta}(x_m)}{\dot{\eta}(x_n)}\widetilde{\mathcal M}_{n\,m},
\quad 
\widetilde{\mathcal M}_{n\,m}\eqdef \frac{\left.
\left(\widetilde{\mathcal H}_{\mathcal{D},\,\mathcal{N}}
\prod_{j\neq m}^{\widetilde{\mathcal N}}\bigl(\eta(x)-\eta_j\bigr)
\right)\right |_{x=x_n}}{\prod_{j\neq n}^{\widetilde{\mathcal N}}\bigl(\eta_n-\eta_j\bigr)},
\label{Mtilde}\\
\mathcal{M}&=-{D}^{-1}\widetilde{\mathcal M}\,{D},\quad
{D}\eqdef \text{diag}(\dot{\eta}(x_1),\ldots,\dot{\eta}(x_{\widetilde{\mathcal N}})).
\label{MMeig}
\end{align}
The eigenvalues of $\mathcal{M}$ are the same as those of $\widetilde{\mathcal M}$ with the sign change
and the corresponding eigenvectors are
\begin{align}   
\widetilde{\mathcal M}\widetilde{\text v}=-\alpha\widetilde{\text v}\quad
\Leftrightarrow  \quad \mathcal{M}\text{v}=\alpha\text{v},\quad
\text{v}\eqdef {D}^{-1}\widetilde{\text v}.
\label{vvrel}
\end{align}

\subsection{Matrix $\widetilde{\mathcal{M}}$ and symmetric matrix $\mathcal{M}$}
\label{sec:mat}
Let us evaluate the matrix elements of $\mathcal{M}$ ($\widetilde{\mathcal M}$) by using the
explicit form of the differential operator $\widetilde{\mathcal H}_{\mathcal D}$ \eqref{ThamD}--\eqref{f2}.
Let us denote the zeros of $\Xi_{\mathcal D}(\eta)$ by $\{\zeta_1,\ldots,\zeta_{\ell_{\mathcal D}}\}$,
which consist of real zeros and possible complex conjugate pairs:
\begin{equation*}
\Xi_{\mathcal D}(\zeta_j)=0,\quad j=1,\ldots,\ell_{\mathcal D},\quad 
\{\zeta_1,\ldots,\zeta_{\ell_{\mathcal D}}\}=\{\zeta_1^*,\ldots,\zeta_{\ell_{\mathcal D}}^*\}\quad
\text{as a set.}
\end{equation*}
By evaluating the equation 
$ \widetilde{\mathcal{H}}_{\mathcal{D},\,\mathcal N} P_{\mathcal{D},\,\mathcal N}(\eta)=0$ 
at one of the zeros $\eta_n$ of $P_{\mathcal{D},\,\mathcal N}(\eta)$,
we obtain the relations among various zeros,
\begin{align}
&2c_2(\eta_n)\sum_{j\neq n}^{\widetilde{\mathcal N}}\frac1{\eta_n-\eta_j}+f_1(\eta_n)=0,\quad
f_1(\eta_n)=c_1(\eta_n,\bm{\lambda}^{[M,N]})-2c_2(\eta_n)\sum_{j=1}^{\ell_{\mathcal D}}
\frac1{\eta_n-\zeta_j},
\label{mulroots}\\
&\qquad \qquad \Rightarrow \sum_{j\neq n}^{\widetilde{\mathcal N}}\frac1{\eta_n-\eta_j}
-\sum_{j=1}^{\ell_{\mathcal D}}
\frac1{\eta_n-\zeta_j}=-\frac{c_1(\eta_n,\bm{\lambda}^{[M,N]})}{2c_2(\eta_n)},\quad 
n=1,\ldots, \widetilde{\mathcal N}.
\label{rrel}
\end{align}
The off-diagonal part of $\widetilde{\mathcal M}$ coming from
$f_1(\eta)\frac{d}{d\eta}$ is
\begin{align*}
&-4f_1(\eta_n)\frac{\left.\sum_{j\neq m}\frac1{\eta-\eta_j}\prod_{k\neq m}(\eta-\eta_k)\right|_{\eta=\eta_n}}
{\prod_{j\neq n}(\eta_n-\eta_j)}
=-4f_1(\eta_n)\frac{\prod_{k\neq n,m}(\eta_n-\eta_k)}{\prod_{j\neq n}(\eta_n-\eta_j)}\n
&\qquad =-4f_1(\eta_n)\frac1{\eta_n-\eta_m}.
\end{align*}
The off-diagonal part of $\widetilde{\mathcal M}$ originating from $c_2(\eta)\frac{d^2}{d\eta^2}$ is
\begin{align*}
(-4)\cdot2c_2(\eta_n)
\frac{\left. \sum_{j\neq n,m}\frac1{\eta-\eta_j}\prod_{l\neq n,m}(\eta-\eta_l)\right|_{\eta=\eta_n}}
{\prod_{j\neq n}(\eta_n-\eta_j)}=(-4)\cdot\frac{2c_2(\eta_n)}{\eta_n-\eta_m}\sum_{j\neq n,m}
\frac1{\eta_n-\eta_j}.
\end{align*}
Summing these two contributions for the off-diagonal part
of $\widetilde{\mathcal M}$ and using \eqref{mulroots}, one obtains
\begin{align}
\widetilde{\mathcal M}_{n\,m}=
\frac{-4}{\eta_n-\eta_m}\left(2c_2(\eta_n)\sum_{j\neq n,m}\frac1{\eta_n-\eta_j}+f_1(\eta_n)\right)
=\frac{8c_2(\eta_n)}{(\eta_n-\eta_m)^2}.
\label{Msym1}
\end{align}
This leads to the  {\em symmetric\/} $\mathcal{M}$ matrix, $\mathcal{M}_{n\,m}=\mathcal{M}_{m\,n}$,
\begin{equation}
\mathcal{M}_{n\,m}=-\frac{\dot{\eta}(x_m)}{\dot{\eta}(x_n)}\frac{8c_2(\eta_n)}{(\eta_n-\eta_m)^2}
=- \frac{2\dot{\eta}(x_n)\dot{\eta}(x_m)}{(\eta_n-\eta_m)^2}=\mathcal{M}_{m\,n},\quad n\neq m,
\label{Msym2}
\end{equation}
since $4c_2(\eta_n)=\bigl(\dot{\eta}(x_n)\bigr)^2$. 
 
Next we evaluate the diagonal  elements of  matrix $\widetilde{\mathcal M}$.
The diagonal elements coming from $f_2(\eta)$ and  the constant are:
\begin{align*}
\widetilde{\mathcal M}_{n\,n}: 
-4\left[ 2c_2(\eta_n)\sum_{j<k}^{\ell_\mathcal{D}}\frac1{\eta_n-\zeta_j}\frac1{\eta_n-\zeta_k}-
c_1(\eta_n,\bm{\lambda}^{[M,N]}-\bm{\delta})\sum_{j=1}^{\ell_{\mathcal D}}\frac1{\eta_n-\zeta_j}
+\frac14\mathcal{E}(\mathcal{N})\right].
\end{align*}
The diagonal part coming from
$f_1(\eta)\frac{d}{d\eta}$ is
\begin{align*}
\widetilde{\mathcal M}_{n\,n}: \quad
&-4f_1(\eta_n)\sum_{k\neq n}^{\widetilde{\mathcal N}}\frac1{\eta_n-\eta_k}\n
=&-4\left(c_1(\eta_n,\bm{\lambda}^{[M,N]})-2c_2(\eta_n)\sum_{j=1}^{\ell_{\mathcal D}}
\frac1{\eta_n-\zeta_j}\right)\sum_{k\neq n}^{\widetilde{\mathcal N}}\frac1{\eta_n-\eta_k}.
\end{align*}
The diagonal part originating from $c_2(\eta)\frac{d^2}{d\eta^2}$ is
\begin{align*}
\widetilde{\mathcal M}_{n\,n}: \quad
(-4)\cdot 2c_2(\eta_n)\sum_{j<k}^{\widetilde{\mathcal N}}{}'\frac1{\eta_n-\eta_j}\frac1{\eta_n-\eta_k},
\end{align*}
in which $\sum'$ means that the singular terms ($j=n$ and $k=n$) are removed.

Summing them, the diagonal element of $\mathcal{M}$ is 
\begin{align}
&\frac14\mathcal{M}_{n\,n}=-\frac14\widetilde{\mathcal M}_{n\,n}\n
&\ 
=c_2(\eta_n)\left(2\sum_{j<k}^{\ell_{\mathcal D}}\frac1{\eta_n-\zeta_j}\frac1{\eta_n-\zeta_k}
-2\sum_{j=1}^{\ell_{\mathcal D}}\frac1{\eta_n-\zeta_j}\sum_{k\neq n}^{\widetilde{\mathcal N}}\frac1{\eta_n-\eta_k}
+2\sum_{j<k}^{\widetilde{\mathcal N}}{}'\frac1{\eta_n-\eta_j}\frac1{\eta_n-\eta_k}\right)\n
&\qquad \qquad -c_1(\eta_n,\bm{\lambda}^{[M,N]}-\bm{\delta})\sum_{j=1}^{\ell_{\mathcal D}}\frac1{\eta_n-\zeta_j}+c_1(\eta_n,\bm{\lambda}^{[M,N]})\sum_{k\neq n}^{\widetilde{\mathcal N}}\frac1{\eta_n-\eta_k}+\frac14\mathcal{E}(\mathcal{N}).
\end{align}
This can be further simplified by using \eqref{rrel}. One factor of 
$-\sum_{j=1}^{\ell_{\mathcal D}}\frac1{\eta_n-\zeta_j}
\sum_{k\neq n}^{\widetilde{\mathcal N}}\frac1{\eta_n-\eta_k}$
 can be changed to 
\begin{align*}
-\sum_{j=1}^{\ell_{\mathcal D}}\frac1{\eta_n-\zeta_j}\left(
\sum_{j=1}^{\ell_{\mathcal D}}\frac1{\eta_n-\zeta_j}-\frac{c_1(\eta_n,\bm{\lambda}^{[M,N]})}{2c_2(\eta_n)}\right),\
\end{align*}
The same factor can be simplified to
\begin{align*}
-\sum_{k\neq n}^{\widetilde{\mathcal N}}\frac1{\eta_n-\eta_k}\left(
\sum_{k\neq n}^{\widetilde{\mathcal N}}\frac1{\eta_n-\eta_k}
+\frac{c_1(\eta_n,\bm{\lambda}^{[M,N]})}{2c_2(\eta_n)}\right).
\end{align*}
We find 
\begin{align*}
&2\sum_{j<k}^{\ell_{\mathcal D}}\frac1{\eta_n-\zeta_j}\frac1{\eta_n-\zeta_k}
-\left(\sum_{j=1}^{\ell_{\mathcal D}}\frac1{\eta_n-\zeta_j}\right)^2
=-\sum_{j=1}^{\ell_{\mathcal D}}\frac1{(\eta_n-\zeta_j)^2},\n
&2\sum_{j<k}^{\widetilde{\mathcal N}}{}'\frac1{\eta_n-\eta_j}\frac1{\eta_n-\eta_k}
-\left(\sum_{k\neq n}^{\widetilde{\mathcal N}}\frac1{\eta_n-\eta_k}\right)^2
=-\sum_{k\neq n}^{\widetilde{\mathcal N}}\frac1{(\eta_n-\eta_k)^2}.
\end{align*}

We finally arrive at
\begin{align}
\frac14\mathcal{M}_{n\,n}&=-\frac14\widetilde{\mathcal M}_{n\,n}\n
&=-c_2(\eta_n)\left(\sum_{j=1}^{\ell_{\mathcal D}}\frac1{(\eta_n-\zeta_j)^2}
+ \sum_{k\neq n}^{\widetilde{\mathcal N}}\frac1{(\eta_n-\eta_k)^2}\right)\n
&\quad+\left(\frac{c_1(\eta_n,\bm{\lambda}^{[M,N]})}{2}-
c_1(\eta_n,\bm{\lambda}^{[M,N]}-\bm{\delta})\right)
\sum_{j=1}^{\ell_{\mathcal D}}\frac1{\eta_n-\zeta_j}\n
&\quad +\frac{c_1(\eta_n,\bm{\lambda}^{[M,N]})}{2}\sum_{k\neq n}^{\widetilde{\mathcal N}}
\frac1{\eta_n-\eta_k}+\frac14\mathcal{E}(\mathcal{N}).
\end{align}
From this we find that 
\begin{align}
\widetilde{\mathcal M}_{n\,n}\in\mathbb{R}\ \text{for real zero}\ \eta_n,\quad
(\widetilde{\mathcal M}_{n\,n})^*=\widetilde{\mathcal M}_{\bar{n}\,\bar{n}} \ 
\text{for complex zero}\ \eta_n,
\end{align}
and for the off-diagonal elements \eqref{Msym1}, \eqref{Msym2} we obtain
\begin{align*}
\widetilde{\mathcal M}_{n\,m}\in\mathbb{R}\ \text{for real zero}\ \eta_n, \eta_m,\quad
&(\widetilde{\mathcal M}_{n\,m})^*=\widetilde{\mathcal M}_{n\,\bar{m}} \ 
\text{for real zero}\ \eta_n, \text{complex zero}\ \eta_m,\\
&(\widetilde{\mathcal M}_{n\,m})^*=\widetilde{\mathcal M}_{\bar{n}\,m} \ 
\text{for complex zero}\ \eta_n, \text{real zero}\ \eta_m,\\
&(\widetilde{\mathcal M}_{n\,m})^*=\widetilde{\mathcal M}_{\bar{n}\,\bar{m}} \
\text{for complex zero}\ \eta_n,  \eta_m.
\end{align*}
With the understanding that $\bar{n}=n$ for  a real zero $\eta_n$, 
the above relations are summarily written as
\begin{equation}
(\widetilde{\mathcal M}_{n\,m})^*=\widetilde{\mathcal M}_{\bar{n}\,\bar{m}}.
\end{equation}

\subsection{Real eigenvalues and orthogonality}
\label{sec:real}
Now we discuss the eigenvalues and eigenvectors of $\widetilde{\mathcal M}$, which act on 
the vector space
$\mathbb{C}^{\widetilde{\mathcal N}}$.
It is easy to see that   $\mathbb{C}^{\widetilde{\mathcal N}}$ is a direct sum 
\begin{equation}
\mathbb{C}^{\widetilde{\mathcal N}}=\widetilde{\bm{V}}\oplus i\widetilde{\bm{V}},\quad
\widetilde{\bm{V}}
\eqdef\{\widetilde{\text v}\in \mathbb{C}^{\widetilde{\mathcal N}}| (\widetilde{\text v}_n)^*
=\widetilde{\text v}_{\bar{n}}\},
\end{equation}
in which $\widetilde{\bm{V}}$ is a vector space over ${\mathbb R}$ and its vector components
are  indexed by the zeros $\{\eta_n\}$ of $P_{\mathcal{D},\,\mathcal{N}}(\eta)$.
Namely $\widetilde{\text v}_n\in\mathbb{ R}$ for a real zero $\eta_n$ and 
$(\widetilde{\text v}_n)^*=\widetilde{\text v}_{\bar{n}}$ for a complex zero $\eta_n$.
Then the vector space $\widetilde{\bm{V}}$ is  invariant under the action of $\widetilde{\mathcal M}$,
\begin{equation}
\widetilde{\mathcal M}\widetilde{\bm{V}}\subseteq \widetilde{\bm{V}}
\Leftrightarrow
(\widetilde{\mathcal M}\widetilde{\text v})_n^*
=(\sum_m \widetilde{\mathcal M}_{n\,m}\widetilde{\text v}_m)^*=
\sum_m\widetilde{\mathcal M}_{\bar{n}\,\bar{m}}\widetilde{\text v}_{\bar{m}}=
\sum_{\bar{m}}\widetilde{\mathcal M}_{\bar{n}\,\bar{m}}\widetilde{\text v}_{\bar{m}}
=(\widetilde{\mathcal M}\widetilde{\text v})_{\bar{n}}.
\end{equation}
This  also means that
all the eigenvalues are real 
\begin{equation}
\widetilde{\mathcal M}\widetilde{\text v}=\alpha \widetilde{\text v}\quad 
\Rightarrow \alpha\in\mathbb{R},
\end{equation}
as $(\widetilde{\mathcal M}\widetilde{\text v})_n$ is real for real zero $\eta_n$.

Next we define another vector space $\bm{V}$ over ${\mathbb R}$, and an inner product 
\begin{equation}
\bm{V}\eqdef\{{\text v}=D^{-1}\widetilde{\text v}|\widetilde{\text v}\in\widetilde{\bm{V}}\},\quad
 (\text{v},\text{w})\eqdef \sum_{m=1}^{\widetilde{\mathcal N}}\text{v}_m\text{w}_m,
 \quad \text{v}, \text{w}\in\bm{V}.
 \end{equation}
It is easy to see that $ (\text{v},\text{w})\in\mathbb{R}$, $\forall\text{v},\forall\text{w}\in\bm{V}$.
The complex symmetry of $\mathcal{M}$ \eqref{Msym2}, instead of the hermiticity (real symmetry)
in the classical case,  leads to the following
\begin{theo}
\label{theo:1}
The  eigenvectors of $\mathcal{M}$ belonging to different eigenvalues are orthogonal to each other
with respect to the above inner product,
\begin{align} 
& 
 \mathcal{M}{\rm v}=\alpha{\rm v},   \quad   
 \mathcal{M}{\rm w}=\beta{\rm w},\quad \alpha\neq \beta, \quad {\rm v}, {\rm w}\in\bm{V},\n
&\alpha({\rm v},{\rm w})=(\mathcal{M}{\rm v},{\rm w})=({\rm v}, \mathcal{M}{\rm w})
=\beta({\rm v},{\rm w})\Rightarrow ({\rm v},{\rm w})=0.
\end{align}
\end{theo}

\bigskip
It is evident that the ordinary lower degree multi-indexed orthogonal polynomials 
$\{P_{\mathcal{D},\, m}(\eta)\}$, $m=0,1,\ldots, \mathcal{N}-1$, provide the 
eigenvectors.
Let us  represent them  
 by Lagrange interpolation formula
\begin{equation*}
P_{\mathcal{D},\, m}\bigl(\eta(x)\bigr)=\sum_{l=1}^{\widetilde{\mathcal N}}
\frac{\prod_{j\neq l}\left(\eta(x)-\eta_j\right)}{\prod_{j\neq l}(\eta_l-\eta_j)}
\cdot P_{\mathcal{D},\, m}(\eta_l),
\quad m=0,\ldots,\mathcal{N}-1,
\end{equation*}
which is exact. By applying $\widetilde{\mathcal H}_{\mathcal{D},\,\mathcal{N}}$ on both sides, 
we obtain
\begin{equation*}
\bigl(\mathcal{E}(m)-\mathcal{E}(\mathcal{N})\bigr)P_{\mathcal{D},\, m}\bigl(\eta(x)\bigr)
=\sum_{l=1}^{\widetilde{\mathcal N}}
\frac{\widetilde{\mathcal H}_{\mathcal{D},\,\mathcal{N}}
\prod_{j\neq l}(\eta(x)-\eta_j)}{\prod_{j\neq l}(\eta_l-\eta_j)}\cdot P_{\mathcal{D},\, m}(\eta_l),
\quad m=0,\ldots,\mathcal{N}-1,
\end{equation*}
which gives rise to the eigenvalue equation for $\widetilde{M}$ by evaluating at $x=x_n$ and dividing 
by $\prod_{j\neq n}(\eta_n-\eta_j)$,
\begin{align} 
\bigl(\mathcal{E}(m)-\mathcal{E}(\mathcal{N})\bigr)
\cdot\frac{P_{\mathcal{D},\, m}(\eta_n)}{\prod_{j\neq n}(\eta_n-\eta_j)}
&=\sum_{l=1}^{\widetilde{\mathcal N}}
\frac{\left.\left(\widetilde{\mathcal H}_{\mathcal{D},\,\mathcal{N}}
\prod_{j\neq l}(\eta(x)-\eta_j)\right)\right|_{x=x_n}}
{\prod_{j\neq n}(\eta_n-\eta_j)}\cdot 
\frac{P_{\mathcal{D},\, m}(\eta_l)}{\prod_{j\neq l}(\eta_l-\eta_j)}, \n
&=\sum_{l=1}^{\widetilde{\mathcal N}}\widetilde{\mathcal M}_{n\, l}\cdot
\frac{P_{\mathcal{D},\, m}(\eta_l)}{\prod_{j\neq l}(\eta_l-\eta_j)}
\quad m=0,\ldots,\mathcal{N}-1.
\end{align}
This leads to the following
\begin{theo}\label{theo:2}
The $\mathcal N$ eigenvalues and eigenvectors of $\mathcal{M}$ corresponding to the  
ordinary lower degree multi-indexed orthogonal polynomials $\{P_{\mathcal{D},\, m}(\eta)\}$ are
\begin{align}
&\mathcal{E}(\mathcal{N})-\mathcal{E}(m),\quad 
\sum_{\ell=1}^{\widetilde{\mathcal N}}\mathcal{M}_{n\,\ell}{\rm v}^{(m)}_{\ell}
=(\mathcal{E}(\mathcal{N})-\mathcal{E}(m))
{\rm v}^{(m)}_n,\quad
m=0,1,\ldots,\mathcal{N}-1.
\label{eigvals}\\
&\qquad\qquad {\rm v}^{(m)}_n\propto \frac{P_{\mathcal{D},\, m}(\eta(x_n))}{\dot{\eta}(x_n)
P'_{\mathcal{D},\,\mathcal N}(\eta(x_n))}
=\frac{P_{\mathcal{D},\, m}(\eta(x_n))}{\quad \left.\left(\frac{d P_{\mathcal{D},\,\mathcal N}(\eta(x))}
{dx}\right)\right|_{x=x_n}},\quad n=1,\ldots, \widetilde{\mathcal N}.
\label{eigveceq}
\end{align}
\end{theo}
This is a proof of \eqref{disor1} in the {\bf Theorem \ref{theo1}}.

\bigskip
It is obvious that the extra lower degree polynomials $\{P_{\mathcal{D}',\,\mathcal{N}}(\eta)\}$
do not satisfy the eigenvalue equations
\begin{equation*}
\widetilde{\mathcal H}_{\mathcal{D}}P_{\mathcal{D}',\,\mathcal{N}}(\eta)
\neq\widetilde{\mathcal E}(\mathcal{D}')\cdot P_{\mathcal{D}',\,\mathcal{N}}(\eta).
\end{equation*}
For the discrete orthogonality, however, it is necessary and sufficient that the above eigenvalue
equations are satisfied only on the zeros
\begin{equation}
\widetilde{\mathcal H}_{\mathcal{D}}P_{\mathcal{D}',\,\mathcal{N}}\left(\eta_j\right)
=\widetilde{\mathcal E}(\mathcal{D}')\cdot P_{\mathcal{D}',\,\mathcal{N}}(\eta_j),
\quad j=1,\ldots, \widetilde{\mathcal N},
\end{equation}
with a certain eigenvalue $\widetilde{\mathcal{\mathcal E}}(\mathcal{D}')$.
This is guaranteed by the following 

\begin{theo}\label{theo:3}
The extra  lower degree polynomials $\{P_{\mathcal{D}',\,\mathcal{N}}(\eta)\}$ satisfy the
equations
\begin{align} 
  \Xi_{\mathcal D}(\eta)\left( 
  \widetilde{\mathcal H}_{\mathcal{D}}P_{\mathcal{D}',\,\mathcal{N}}(\eta)
-\widetilde{\mathcal E}(\mathcal{D}')\cdot
P_{\mathcal{D}',\,\mathcal{N}}(\eta)\right)
=\mathcal{E}'(\mathcal{D}')P_{\mathcal{D},\,\mathcal{N}}(\eta)
\Xi_{\mathcal{D}'}(\eta),
\label{basic}
\end{align}
in which the r.h.s. are proportional to the product of the
multi-indexed polynomial $P_{\mathcal{D},\,\mathcal{N}}(\eta)$
and the denominator polynomial $\Xi_{\mathcal{D}'}(\eta)$ corresponding to the multi-index
of the extra lower degree polynomial $P_{\mathcal{D}',\,\mathcal{N}}(\eta)$.
The eigenvalues $\{\widetilde{\mathcal{\mathcal E}}(\mathcal{D}')\}$ and the coefficients
$\{\mathcal{E}'(\mathcal{D}')\}$ 
are listed in Theorem {\bf\rm \ref{EDrelation}}.
\end{theo}

The proof is quite elementary. As remarked at the end  of section 2, the second order differential operator 
$\Xi_{\mathcal{D}} \widetilde{\mathcal{H}}_{\mathcal{D}}$ maps any polynomial (in $\eta$)
to  another polynomial. Since the multi-index polynomials and the denominator polynomials are given explicitly 
for any choices of $\mathcal{D}$ and $\mathcal{D}'$, it is straightforward to verify the equality of (4.30) by
evaluating both sides explicitly.

The presence of $\Xi_{\mathcal{D}'}(\eta)$ on r.h.s. is easily guessed by degree counting.
By evaluating the above basic equations at the zeros and following the same logic as that for the ordinary
lower degree polynomials, we obtain the eigenvectors of $\widetilde{\mathcal M}$
corresponding to the extra lower degree polynomials $\{P_{\mathcal{D}',\,\mathcal{N}}(\eta)\}$:
\begin{align} 
&\bigl(\widetilde{\mathcal E}(\mathcal{D}')-\mathcal{E}(\mathcal{N})\bigr)
\cdot\frac{P_{\mathcal{D}',\, \mathcal{N}}(\eta_n)}{\prod_{j\neq n}(\eta_n-\eta_j)}
=\sum_{l=1}^{\widetilde{\mathcal N}}\widetilde{\mathcal M}_{n\, l}\cdot
\frac{P_{\mathcal{D}',\,  \mathcal{N}}(\eta_l)}{\prod_{j\neq l}(\eta_l-\eta_j)},\n
 &{\widetilde{\text v}}^{\mathcal{D}'}_n\eqdef
 \frac{P_{\mathcal{D}',\, \mathcal{N}}(\eta_n)}{\prod_{j\neq n}(\eta_n-\eta_j)}
\ \Rightarrow \ 
{\rm v}^{\mathcal{D}'}_n\propto \frac{P_{\mathcal{D}',\, \mathcal{N}}(\eta(x_n))}{\dot{\eta}(x_n)
P'_{\mathcal{D},\,\mathcal N}(\eta(x_n))}
=\frac{P_{\mathcal{D}',\, \mathcal{N}}(\eta(x_n))}{\quad \left.\left(\frac{d P_{\mathcal{D},\,\mathcal N}(\eta(x))}
{dx}\right)\right|_{x=x_n}}.
\label{eigveceq2}
\end{align}
This is a proof of \eqref{disor2} and \eqref{disor3} in the {\bf Theorem \ref{theo1}}.
\begin{coro}
\label{chris4}
The above eigenvector formulas \eqref{eigveceq},\eqref{eigveceq2}  show the forms of the Christoffel
numbers
\begin{equation}
\lambda_n\propto \frac1{\dot{\eta}(x_n)^2\bigl(P_{\mathcal{D},\,\mathcal{N}}'(\eta_n)\bigr)^2},\quad
(\lambda_n)^*=\lambda_{\bar{n}}.
\end{equation}
\end{coro}
This means that Gauss quadrature formula can be used for real functions by using the
zeros of $P_{\mathcal{D},\,\mathcal{N}}(\eta)$. 
It would be interesting to find out characterisation of functions for which using multi-indexed polynomials
is more effective than other orthogonal polynomials.

The eigenvalues $\widetilde{\mathcal E}(\mathcal{D}')$ depend on the deleted ($d_j^{\ai,\ait}$) and
added ($\epsilon_k^{\ai,\ait}$) degrees. They are neatly expressed by the complements of 
$\mathcal{D}'$ in $\mathcal{D}$ and vice versa, as shown by the following

\begin{theo}
\label{EDrelation}
For the three types of extra `lower degree polynomials' \eqref{1D}, \eqref{2D}, \eqref{3D}, 
the corresponding eigenvalues of $\widetilde{\mathcal M}$ are
\begin{equation*}
\widetilde{\mathcal E}(\mathcal{D}'_{{\lambda}, j,k})-\mathcal{E}(\mathcal{N}),\quad \mathcal{D}'\subseteq\mathcal{EP},~~
\lambda=\ai,\ait,\at,
\end{equation*}
and
\begin{align}
\widetilde{\mathcal E}\bigl(\mathcal{D}'_{\ai,j,k}\bigr)&\eqdef
\widetilde{\mathcal E}^\ai\bigl(d^\ai_j\bigr) +\widetilde{\mathcal E}^{\ai}\bigl(\epsilon^\ai_k\bigr)
 -\mathcal{E}(\mathcal{N}), \quad \quad \, \rm{type\ I},
\label{eigvvec1}
\\
\widetilde{\mathcal E}\bigl(\mathcal{D}'_{\ait,j,k}\bigr)&\eqdef
\widetilde{\mathcal E}^\ait\bigl(d^\ait_j\bigr) +\widetilde{\mathcal E}^{\ait}\bigl(\epsilon^\ait_k\bigr)
-\mathcal{E}(\mathcal{N}), \quad \rm{type\ II},
\label{eigvvec2}
\\
\widetilde{\mathcal E}\bigl(\mathcal{D}'_{\at,j,k}\bigr)&\eqdef
\widetilde{\mathcal E}^\ai\bigl(d^\ai_j\bigr) +\widetilde{\mathcal E}^{\ait}\bigl(d^\ait_k\bigr)
-\mathcal{E}(\mathcal{N}), \quad \rm{type\ III},
\label{eigvvec3}
\end{align}
in which the energies of the virtual state solutions $\widetilde{\mathcal E}^\ai(\rm{v})$,
$\widetilde{\mathcal E}^\ait(\rm{v})$ are given in \S\ref{sec:ori} for \text{\rm L} and \text{\rm J}.
The coefficients $\mathcal{E}'(\mathcal{D}')$  on r.h.s. of the basic equations \eqref{basic} are:
\begin{align}
& 
2\left(\mathcal{E}(\mathcal{N})-\widetilde{\mathcal E}^\ai(\epsilon^\ai_k)\right),
\quad \ \rm{L}, \rm{J} \qquad \qquad \qquad \quad \ \  \rm{type\ I},
\\
&
 2\left(\mathcal{E}(\mathcal{N})-\widetilde{\mathcal E}^\ait(\epsilon^\ait_k)\right),
\quad \rm{L}, \rm{J} \qquad \qquad \qquad  \quad \ \rm{type\ II},
\\
&
-8\ \rm{for\ L}, \ -32\ \rm{for\ J},\ \rm{or} -2 c_\mathcal{F}^2\ \quad \rm{for\ L,\ J }\ \rm{type\ III}.
\end{align}
\end{theo}


In order to bring out clearly the symmetry in the expressions (\ref{eigvvec1})-(\ref{eigvvec3}), hereafter we shall ``abuse" the notations somewhat.
In set theory, $\mathcal{D}^I\backslash\mathcal{D}^{\prime I}$ means a set $\{d^I_j\}$, and
$\mathcal{D}^{\prime I}\backslash\mathcal{D}^I$ means a set $\{\epsilon^I_j\}$, etc.
From now on we shall understand  $\mathcal{D}^I\backslash\mathcal{D}^{\prime I}$ as an argument of 
$\widetilde{\mathcal E}$ meaning the element $d^I_j$, and $\mathcal{D}^{\prime I}\backslash\mathcal{D}^I$ meaning  the element  $\epsilon^I_j$. With this understanding, eqs.\, (\ref{eigvvec1})-(\ref{eigvvec3}) become
\begin{align}
\widetilde{\mathcal E}\bigl(\mathcal{D}'_{\ai,j,k}\bigr)&=
\widetilde{\mathcal E}^\ai\bigl(\mathcal{D}^{\ai}\backslash\mathcal{D}^{\prime\ai}\bigr)
 +\widetilde{\mathcal E}^{\ai}\bigl(\mathcal{D}^{\prime\ai}\backslash\mathcal{D}^{\ai}\bigr)
 -\mathcal{E}(\mathcal{N}), \quad \qquad \rm{type\ I},
\label{eigvvec1a}
\\
\widetilde{\mathcal E}\bigl(\mathcal{D}'_{\ait,j,k}\bigr)&=
\widetilde{\mathcal E}^\ait\bigl(\mathcal{D}^{\ait}\backslash\mathcal{D}^{\prime\ait}\bigr)
+\widetilde{\mathcal E}^\ait\bigl(\mathcal{D}^{\prime\ait}\backslash\mathcal{D}^{\ait}\bigr)
-\mathcal{E}(\mathcal{N}), \quad \rm{type\ II},
\label{eigvvec2a}
\\
\widetilde{\mathcal E}\bigl(\mathcal{D}'_{\at,j,k}\bigr)&=
\widetilde{\mathcal E}^\ai\bigl(\mathcal{D}^{\ai}\backslash\mathcal{D}^{\prime\ai}\bigr)
+\widetilde{\mathcal E}^\ait\bigl(\mathcal{D}^{\ait}\backslash\mathcal{D}^{\prime\ait}\bigr)
-\mathcal{E}(\mathcal{N}), \quad \ \,\rm{type\ III},
\label{eigvvec3a}
\end{align}


The orthogonality of the eigenvectors of $\mathcal{M}$, except for the cases of degeneracy,
\begin{equation}
(\text{v}^{\mathcal{D}'},\text{v}^{\mathcal{D}''})=(\text{v}^{\mathcal{D}'},\text{v}^{(m)})
=(\text{v}^{(n)},\text{v}^{(m)})=0,\quad \mathcal{D}'\neq\mathcal{D}'',\quad n\neq m,
\end{equation}
is also verified directly,  independently of the discrete orthogonality relations 
\eqref{disor1}--\eqref{disor3}.
The explicit forms of the eigenvalues are essential to identify the extra `lower degree polynomials',
or the corresponding multi-indices $\mathcal{D}'\subseteq\mathcal{EP}$. 
The structure of the extra `lower degree polynomials' $\mathcal{D}'_{\ai}$,
$\mathcal{D}'_{\ait}$, $\mathcal{D}'_{\at}$, \eqref{eigvvec1}--\eqref{eigvvec3}, 
 clarify the meaning of the formula 
 \eqref{ellD} for  the extra dimensions $\ell_{\mathcal D}$.

\bigskip
It is important to stress that the basic equations \eqref{basic} are polynomials in the variable $\eta$ 
and in the parameters $g$, $h$.
The equations involve  differentiation, multiplications and additions of the polynomials only.

\bigskip
It is interesting to note the symmetry in $\mathcal{D}$ and $\mathcal{D}'$  of the eigenvalues 
\eqref{eigvvec1a}--\eqref{eigvvec3a} 
in the basic equations \eqref{basic}.
In fact, the basic equations \eqref{basic} hold with slight modifications when 
$\mathcal{D}$ and $\mathcal{D}'$ are exchanged as shown in the following 
\begin{theo}\label{theo:6}
For a given multi-index $\mathcal{D}$ \eqref{defD}, we have a set of multi-indices 
$\{\mathcal{D}'\}$ corresponding to the extra lower degree polynomials \eqref{1D}--\eqref{3D}.
For each multi-index $\mathcal{D}'$ belonging to $\mathcal{D}$ the following equation holds
\begin{align} 
  \Xi_{\mathcal{D}'}(\eta)\left( 
  \widetilde{\mathcal H}_{\mathcal{D}'}P_{\mathcal{D},\,\mathcal{N}}(\eta)
-\widetilde{\mathcal E}(\mathcal{D})\cdot
P_{\mathcal{D},\,\mathcal{N}}(\eta)\right)
=\mathcal{E}^{\prime\prime}(\mathcal{D})P_{\mathcal{D}',\,\mathcal{N}}(\eta)
\Xi_{\mathcal{D}}(\eta),
\label{rbasic}
\end{align}
with the same eigenvalues $\widetilde{\mathcal E}(\mathcal{D})=\widetilde{\mathcal E}(\mathcal{D}')$ 
as \eqref{basic} but with modified coefficients $\mathcal{E}^{\prime\prime}(\mathcal{D})$ 
on the r.h.s. given by
\begin{align}
2\left(\mathcal{E}(\mathcal{N})-\widetilde{\mathcal E}^\ai(\mathcal{D}\backslash\mathcal{D}'\right),
\quad \rm{L}, \rm{J} \qquad \ \,  &\rm{type\ I},\\
 2\left(\mathcal{E}(\mathcal{N})-\widetilde{\mathcal E}^\ait(\mathcal{D}\backslash\mathcal{D}')\right),
\quad \rm{L}, \rm{J} \qquad\ &\rm{type\ II},\\
-\frac{1}{\gamma}
\left(\mathcal{E}(\mathcal{N})-\widetilde{\mathcal E}^\ai(\mathcal{D}^\ai\backslash\mathcal{D}^{'\ai})\right)\!\!
\left(\mathcal{E}(\mathcal{N})-\widetilde{\mathcal E}^\ait(\mathcal{D}^\ait\backslash\mathcal{D}^{'\ait})\right),
~~ \gamma= \frac{c_\mathcal{F}^2}{2},
\quad \rm{L}, \rm{J} \quad &\rm{type\ III}.
\end{align}
\end{theo}
These eigenvalue formulas are obtained by educated guess based on the explicit evaluations of the eigenvalues of 
$\widetilde{\mathcal M}$.

It is interesting to note that formal replacement $\mathcal{D}'\to\mathcal{D}$,
with $\mathcal{D}'\backslash\mathcal{D}=\phi=\mathcal{D}\backslash\mathcal{D}'$,
$\widetilde{\mathcal E}^{\ai,\ait}(\phi)=0$, reduces the basic equations \eqref{basic}, \eqref{rbasic} 
to the ordinary eigenvalue equation
$\widetilde{\mathcal H}_{\mathcal{D}}P_{\mathcal{D},\,\mathcal{N}}(\eta)=
\mathcal{E}(\mathcal{N})P_{\mathcal{D},\,\mathcal{N}}(\eta)$ \eqref{fuchs} for the type I and II cases.
These theorems show rich and interesting structures of multi-indexed
Laguerre and Jacobi polynomials.

\section{Krein-Adler orthogonal polynomials}
\label{sec:KA}

In this section we demonstrate the discrete orthogonality relations for another family of
`new' orthogonal polynomials, that is, the Krein-Adler polynomials based on the
Hermite, Laguerre and Jacobi polynomials.
They have simpler structures and longer histories than the exceptional/multi-indexed
polynomials. Most of the Theorems in section \ref{sec:disortmulLJ} and \ref{sec:zerosum} 
hold with minor modifications.
\subsection{Overview of the Krein-Adler  polynomials}
\label{sec:KAover}

The K-A systems based on the Hermite, Laguerre and Jacobi
systems are the iso-spectral deformations of the original systems \S\ref{sec:oriher}--\S\ref{sec:orijac} through
multiple Darboux transformations using the eigenfunctions specified by the multi-index
$\mathcal{D}$ \cite{krein,adler}:

\begin{align}
&\mathcal{D}=\{d_1,\ldots,d_M\},\qquad
0<d_1<\cdots<d_M,\n
  U_{\mathcal{D}}(x)&\eqdef U(x)-2\partial_x^2
  \log\bigl|\text{W}[\phi_{d_1},\ldots,\phi_{d_M}](x)\bigr|,
  \label{KApot}\\
  \phi_{\mathcal{D},\,n}(x)&\eqdef
  \frac{\text{W}[\phi_{d_1},\ldots,\phi_{d_M},\phi_n](x)}
  {\text{W}[\phi_{d_1},\ldots,\phi_{d_M}](x)}
  \hspace{40mm} \ (n\in\mathbb{Z}_{\ge0}\backslash\mathcal{D})\\
 &=\phi_0(x){\eta'(x)}^M\frac{\text{W}[P_{d_1},\ldots,P_{d_M},P_n](\eta)}
  {\text{W}[P_{d_1},\ldots,P_{d_M}](\eta)}
=\phi_0(x){\eta'(x)}^M\!\!\!\cdot\!P_{\mathcal{D},n}(\eta)/{\Xi_{\mathcal D}(\eta)},
\label{KAeig}\\
 P_{\mathcal{D},n}(\eta)&\eqdef \text{W}[P_{d_1},\ldots,P_{d_M},P_n](\eta),\quad
 \Xi_{\mathcal D}(\eta)\eqdef \text{W}[P_{d_1},\ldots,P_{d_M}](\eta),
 \label{KApols}\\
  \phi_{\mathcal{D},\,n}(x) &=\psi_{\mathcal{D}}(x)
  P_{\mathcal{D},n}\bigl(\eta(x)\bigr),\quad
  \psi_{\mathcal{D}}(x)\eqdef\phi_0(x){\eta'(x)}^M/
  {\Xi_{\mathcal{D}}\bigl(\eta(x)\bigr)},\\
\mathcal{H}_{\mathcal D} \phi_{\mathcal{D},\,n}(x)&
=\mathcal{E}(n) \phi_{\mathcal{D},\,n}(x),
\hspace{58mm} \ (n\in\mathbb{Z}_{\ge0}\backslash\mathcal{D}).
\label{eigproKA}
\end{align}
In this case there is no distinction of type I and II.
For the non-singularity of the potential $U_{\mathcal D}(x)$ \eqref{KApot} and the positivity
of the norms of the eigenfunctions \eqref{KAeig}, the multi-index $\mathcal{D}$ is required to
satisfy the Krein-Adler conditions \cite{krein,adler},
\begin{equation}
\prod_{j=1}^M(m-d_j)\ge0, \quad \forall m\in\mathbb{Z}_{\ge0}.
\label{KAcond}
\end{equation}
Here again, although we have restricted $\{d_j\}$ to be positive, including the case of 
$d_1=0$ goes almost parallel.

Since all the polynomial entries for the Wronskians in 
$P_{\mathcal{D},n}(\eta)$ and $ \Xi_{\mathcal D}(\eta)$
\eqref{KApols} are the original ones $\{P_n(\eta)\}$, various calculations become
much simpler than those for the multi-indexed cases.
For generic values of the  parameters $g$ and $h$, the K-A polynomial 
$P_{\mathcal{D},n}(\eta)$ is of degree $\ell'_{\mathcal D}+n$
and the denominator polynomial $\Xi_{\mathcal D}(\eta)$ is of degree $\ell'_{\mathcal D}+M$ in
$\eta$, in which $\ell'_{\mathcal D}$ is given by
\begin{equation}
\ell_{\mathcal D}'\eqdef \sum_{j=1}^Md_j -\frac12 M(M+1).
 \label{KAellD}
\end{equation}
 
 The second order differential operator
$\widetilde{\mathcal{H}}_{\mathcal{D}}$ obtained from $\mathcal{H}_{\mathcal D}$ 
by a similarity transformation in terms of $\psi_{\mathcal D}(x)$
\begin{equation*}
  \widetilde{\mathcal{H}}_{\mathcal{D}} \eqdef\psi_{\mathcal{D}}(x)^{-1}\circ
  \mathcal{H}_{\mathcal{D}}\circ
  \psi_{\mathcal{D}}(x),\qquad
  \widetilde{\mathcal{H}}_{\mathcal{D}}
  P_{\mathcal{D},n}(\eta)=\mathcal{E}(n)
  P_{\mathcal{D},n}(\eta),
\end{equation*}
has essentially the same form as that for the multi-indexed polynomials \eqref{ThamD}--\eqref{parashift}.
For L and J, the functions $c_1$ and $c_2$ are the same as given in \eqref{cF,c1,c2} with $N=0$.
For the Hermite family these functions are,
\begin{equation}
\text{H}: \quad c_1(\eta)=-\eta/2,\quad c_2=1/4.
\end{equation}
The only difference is that function $f_2(\eta)$ \eqref{f2} has an additional constant term
\begin{equation}
\text{H}: -M/2, \quad \text{L}: -M, \quad \text{J}: -M(g+h+M),
\quad \text{or} \ -\frac14 \mathcal{E}(M) ~~\ \text{collectively}.
\end{equation}

\subsection{Discrete orthogonality relations }
\label{sec:KAdisc}

For a given positive integer $\mathcal{N}\notin\mathcal{D}$,  we formulate
the discrete orthogonal relations in terms of the zeros of 
$P_{\mathcal{D},\,\mathcal{N}}(\eta)$.
\begin{equation*}
P_{\mathcal{D},\,\mathcal{N}}(\eta_j)=0,\quad 
j=1,\ldots,\widetilde{\mathcal N}\eqdef \ell'_{\mathcal D}+\mathcal{N}.
\end{equation*}

The ordinary lower degree polynomials are
\begin{equation}
\{P_{\mathcal{D},\,m}(\eta)\},\quad
m\in\{0,1,\ldots,\mathcal{N}-1\}\backslash\mathcal{D}.
\label{KAord}
\end{equation}
There are $\mathcal{N}-K$ such polynomials with $K\eqdef \#\{d_j\in\mathcal{D}|d_j<\mathcal{N}\}$.
Due to the Wronskian structure, the extra lower degree polynomials have essentially the same structure 
as before, except that the replaced term $\epsilon_k$ should not coincide with $\mathcal{N}$, that is,
\begin{equation}
\{P_{\mathcal{D}',\,\mathcal{N}}(\eta)\},\qquad
\mathcal{D}'\eqdef 
\{d_1,\ldots,d_{j-1},\epsilon_k,d_{j+1},\ldots,d_M\},
\quad d_j\in \mathcal {D},
\quad \mathcal{N}\neq\epsilon_k \in {E}_j.
\label{KA1D}
\end{equation}
It is easy to see that the number of the extra lower degree polynomials is $\ell_{\mathcal D}'+K$, resulting in  the  total number of the lower degree polynomials $\ell_{\mathcal D}'+\mathcal{N}$. 
They satisfy the  discrete orthogonality relations presented in {\bf Theorem \ref{theo1}}
 except that
the degrees $m$ and $n$ in \eqref{disor1}, \eqref{disor2} should read
$m\neq n\in\{0,1,\ldots,\mathcal{N}-1\}\backslash\mathcal{D}$.

The proof of the discrete orthogonality goes parallel with that in Sect. 4.
The same structures of the second order differential operators 
$\widetilde{\mathcal H}_{\mathcal D}$ means that the properties of 
the corresponding matrices $\mathcal{M}$ and $\widetilde{\mathcal M}$ 
\eqref{Mtilde}, \eqref{MMeig} are the same.
That is, {\bf Theorem \ref{theo:1}} holds and the orthogonality of the
eigenvectors is guaranteed. {\bf Theorem \ref{theo:2}} shows that the ordinary lower degree
polynomials \eqref{KAord}  form eigenvectors of $\mathcal{M}$ with the eigenvalue 
$\mathcal{E}(\mathcal{N})- \mathcal{E}(m)$.
The extra lower degree
polynomials \eqref{KA1D}  also provide the eigenvectors of $\mathcal{M}$ with the eigenvalues 
\begin{equation}
{\mathcal E}(\mathcal{N})- \widetilde{\mathcal E}(\mathcal{D}'),\qquad
\widetilde{\mathcal E}\bigl(\mathcal{D}'\bigr)\eqdef
{\mathcal E}\bigl(\mathcal{D}\backslash\mathcal{D}'\bigr)
+{\mathcal E}\bigl(\mathcal{D}'\backslash\mathcal{D}\bigr)
-\mathcal{E}(\mathcal{N}),
\label{KAeig2}
\end{equation}
which looks similar to \eqref{eigvvec1a} and \eqref{eigvvec2a}.
The elements of the eigenvalue formulas are all the eigenvalues $\{\mathcal{E}(n)\}$ 
of the original systems, H, L, J, \S\ref{sec:oriher}--\S\ref{sec:orijac},  instead of the energies 
of the virtual state solutions $\{\widetilde{\mathcal E}^{\ai,\ait}(n)\}$ 
appearing in \eqref{eigvvec1},\eqref{eigvvec2}.
The extra  lower degree polynomials $\{P_{\mathcal{D}',\,\mathcal{N}}(\eta)\}$ satisfy the
basic equations
\begin{align} 
  \Xi_{\mathcal D}(\eta)\left( 
  \widetilde{\mathcal H}_{\mathcal{D}}P_{\mathcal{D}',\,\mathcal{N}}(\eta)
-\widetilde{\mathcal E}(\mathcal{D}')\cdot
P_{\mathcal{D}',\,\mathcal{N}}(\eta)\right)
=\mathcal{E}'(\mathcal{D}')P_{\mathcal{D},\,\mathcal{N}}(\eta)
\Xi_{\mathcal{D}'}(\eta),
\label{KAbasic}
\end{align}
with the coefficients
\begin{equation}
\mathcal{E}'(\mathcal{D}')\eqdef2\bigl(\mathcal{E}(\mathcal{N})-
\mathcal{E}(\mathcal{D}'\backslash\mathcal{D})\bigr).
\end{equation}
The basic equations with $\mathcal{D}$ and $\mathcal{D}'$ exchanged also hold
\begin{align} 
  \Xi_{\mathcal{D}'}(\eta)\left( 
  \widetilde{\mathcal H}_{\mathcal{D}'}P_{\mathcal{D},\,\mathcal{N}}(\eta)
-\widetilde{\mathcal E}(\mathcal{D})\cdot
P_{\mathcal{D},\,\mathcal{N}}(\eta)\right)
=\mathcal{E}^{\prime\prime}(\mathcal{D})P_{\mathcal{D}',\,\mathcal{N}}(\eta)
\Xi_{\mathcal{D}}(\eta),
\label{KArbasic}
\end{align}
with the same eigenvalues $\widetilde{\mathcal E}(\mathcal{D})=\widetilde{\mathcal E}(\mathcal{D}')$, 
\eqref{KAeig2} but with  slightly modified coefficients $\mathcal{E}^{\prime\prime}(\mathcal{D})$
\begin{equation}
\mathcal{E}^{\prime\prime}(\mathcal{D})\eqdef2\bigl(\mathcal{E}(\mathcal{N})-
\mathcal{E}(\mathcal{D}\backslash\mathcal{D}')\bigr).
\end{equation}

The proofs of these basic equations (\ref{KAbasic}) and (\ref{KArbasic})  are again elementary calculation of polynomials.

\begin{coro}
\label{dioph}
Diophantine properties. The matrices $\mathcal{M}$ and $\widetilde{\mathcal M}$ for the
K-A Hermite and Laguerre systems have remarkable diophantine properties. 
Their eigenvalues are all integers, irrespective of the values of the parameter $g$ in the K-A L systems,
as obvious from \eqref{KAeig2}.
For integer and/or half integer parameters $g$ and $h$, the eigenvalues of the  
matrices $\mathcal{M}$ and $\widetilde{\mathcal M}$ for the K-A Jacobi and the multi-indexed 
Laguerre and Jacobi systems are all integers, see \eqref{eigvvec1}--\eqref{eigvvec3}, \eqref{KAeig2}.
\end{coro}

This concludes the proof of the discrete orthogonality for the Krein-Adler polynomials based
on the Hermite, Laguerre and Jacobi polynomials.
\section{Summary and Comments}
\label{sec:sum}

In this paper we have demonstrated the rich and interesting structures of the multi-indexed Laguerre 
and Jacobi polynomials and the K-A polynomials through the discrete orthogonality relations. 
These new orthogonal polynomials $P_{\mathcal{D},\,\mathcal{N}}(\eta)$ 
have extra zeros outside the orthogonality domain.
By expanding the second order differential operator 
$\widetilde{\mathcal H}_{\mathcal D}$ \eqref{ThamD} governing these 
new orthogonal polynomials around their zeros, two 
$\widetilde{\mathcal N}\times \widetilde{\mathcal N}$ matrices ${\mathcal M}$ and
$\widetilde{\mathcal M}$ \eqref{Mtilde}, \eqref{MMeig} are obtained.
Here $\widetilde{\mathcal N}$ is the number of zeros, the ordinary plus the extra.
The complex symmetry of ${\mathcal M}$ leads to the orthogonality of its eigenvectors,
{\bf Theorem \ref{theo:1}}.
The ordinary lower degree polynomials, being the eigenfunctions of 
$\widetilde{\mathcal H}_{\mathcal D}$,
naturally provide the eigenvectors of ${\mathcal M}$ through Lagrange interpolation formula.
The non-eigenfunction parts of the extra lower degree polynomials are shown to be proportional
to the starting function $P_{\mathcal{D},\,\mathcal{N}}(\eta)$ through the basic equations,
{\bf Theorem \ref{theo:3}}.
This provides the eigenvectors of ${\mathcal M}$ corresponding to the extra lower degree polynomials.
The compositions of the extra lower degree polynomials \eqref{1D}--\eqref{3D} and 
{\bf Theorem \ref{theo1}} are quite natural due to the Wronskian construction of these
new orthogonal polynomials  \eqref{PDn}, \eqref{KApols}.
It should be stressed that essentially the same formulas apply to the multi-indexed Laguerre and Jacobi 
polynomials as well as to the K-A polynomials based on the Hermite, Laguerre and Jacobi polynomials.

For mathematical logics only, one could have started from the eigenvalue problem of matrix 
$\mathcal{M}$ ($\widetilde{\mathcal{M}}$) \eqref{Mtilde}-\eqref{MMeig},
without the subtlety of the assumption of infinitesimal oscillation around the zeros of 
$P_{\mathcal{D},\,\mathcal{N}}\bigl(\eta(x)\bigr)$ \eqref{ansatz}.
Then the presence of $\dot{\eta}(x)$ in \eqref{Mtilde} needs `explanation'.
Without $\dot{\eta}(x_n)$ the symmetry of $\mathcal{M}$ is not achieved and 
 the orthogonality does not hold.
As stressed in \cite{os7,os12,os13,os24}, the `sinusoidal coordinate' $\eta(x)$ is essential for the
unified treatment of the classical and multi-indexed orthogonal polynomials,
which could be termed as the `quantum mechanical formulation of the classical orthogonal polynomials'.
The main part of this process could be considered as blow-ups of the second order differential operator
$\widetilde{\mathcal H}_{\mathcal D}$ \eqref{ThamD}, 
around the zeros of these polynomials at degree $\mathcal{N}$ ($\widetilde{\mathcal N}$).

A slightly different formulation of the matrix $\widetilde{\mathcal{M}}$ and its eigenvalue problem, 
which was applicable to non-classical orthogonal polynomials, was discussed in \cite{bihun}.
Interesting  applications to Krall-Laguerre and Krall-Jacobi polynomials were reported. 

The set of extra lower degree polynomials $\{P_{\mathcal{D}',\,\mathcal{N}}(\eta)\}$ for 
$\mathcal{N}=0$ would correspond to the
`exceptional invariant subspace' \cite{xsp}. 
It would be interesting to explore the geometric properties of this subspace
together with the $\ell_{\mathcal D}$ zeros of the denominator polynomial 
$\Xi_{\mathcal D}(\eta)$\eqref{XiD}, 
which are the apparent singularities of the second order differential operator 
$\widetilde{\mathcal H}_{\mathcal D}$ \eqref{ThamD}.

For the actual verifications of the discrete orthogonality {\bf Theorem \ref{theo1}} by 
numerical evaluations, 
 the bounds \eqref{Lbound},  \eqref{Jbound}
on the  parameters $g$ and $h$ or the Krein-Adler conditions on the multi-indices \eqref{KAcond}
do not constitute any severe limitations, as the orthogonality is the consequences of the fact
that the polynomials satisfy the second order basic equations 
\eqref{basic}, \eqref{rbasic}, \eqref{KAbasic}, \eqref{KArbasic}. 
It is, however,  necessary to avoid the exceptional values of the  parameters $g$ and $h$ for which
the highest degree component of $P_{\mathcal{D},\,n}(\eta)$ vanishes.

It is expected that the same theorems would apply to the multi-indexed and/or the K-A versions 
of the Racah, $q$-Racah, 
Wilson and Askey-Wilson polynomials \cite{dismul} etc.
It is a good challenge to demonstrate the discrete orthogonality relations for these 
multi-indexed polynomials.

We have 
realised that the discrete orthogonality
relations for the multi-indexed and K-A systems also hold for any complex values 
of the parameters $g$ and $h$.
This is a simple consequence of the complex {\em symmetry\/} of the matrix $\mathcal{M}$, \eqref{Msym2},
which is derived without any restrictions on the parameters.
The basic equations \eqref{basic}, \eqref{rbasic}, \eqref{KAbasic}, \eqref{KArbasic}, are polynomials in the
parameters $g$ and $h$. They guarantee to make the  lower degree polynomials as the eigenvectors 
of $\mathcal{M}$ for any complex values of the parameters.
Let us emphasise that the quantum mechanical formulation is essential for these results,
as the existence of $\dot{\eta}(x)$ makes $\mathcal{M}$ symmetric in \eqref{Msym2}.
In terms of polynomials only, one goes as far as $\widetilde{\mathcal M}$ \eqref{Mtilde}.

\section*{Acknowledgements}
R.\,S. thanks P-M. Ho for discussion and hospitality at Dept. Physics, National Taiwan University,
where a part of this work was done. 
R.\,S. also  wishes to thank S.-H. Lai, J.-C. Lee and S. Odake for discussions in the early stage of this work.
This work is supported in part by the Ministry of Science and Technology (MoST)
of the Republic of China under Grant  MOST 107-2112-M-032-002 and MOST 109-2112-M-032-008.

\bigskip
\centerline{\bf Availability of data}


$\bullet$ The data that support the  findings of this study are available from the corresponding author upon reasonable request.

%


\end{document}